\documentclass[11pt,french]{amsart}
\usepackage{latexsym}
\usepackage[all]{xy}
\usepackage{amsmath}
\usepackage{amssymb}
\usepackage{amsfonts}

\newcommand  {\colim}  {\mathop{\mathrm{colim}}}

\newcommand  {\St}   {\mathbf{St}}

\newcommand{\s}{\infty}

\newtheorem{thm}{Th\'eor\`eme}[section]
\newtheorem{prop}[thm]{Proposition}
\newtheorem{lem}[thm]{Lemme}

\newtheorem{df}[thm]{D\'efinition}
\newtheorem{cor}[thm]{Corollaire}

\begin{document}

\title{\textbf{Syst\`emes de points dans les dg-cat\'egories satur\'ees}}
\bigskip
\bigskip

\author{Bertrand To\"en}
   \address{IMT UMR 5219, CNRS, Universit\'{e} Paul Sabatier, Toulouse - France}
   \email{bertrand.toen@math.univ-toulouse.fr}
\author{Michel Vaqui\'e}
   \address{IMT UMR 5219, CNRS, Universit\'{e} Paul Sabatier, Toulouse - France}
   \email{michel.vaquie@math.univ-toulouse.fr}

\bigskip 

\date{Mars 2015}

\dedicatory{\`A Vadim Schechtman, avec admiration et amiti\'e.}

\maketitle

\begin{abstract}
Dans ce travail nous consid\'erons le probl\`eme de r\'ealiser
g\'eom\'etriquement les cat\'egories triangul\'ees (plut\^ot les dg-cat\'egories
triangul\'ees) comme des cat\'egories d\'eriv\'ees de vari\'et\'es alg\'ebriques.
Pour cela, on introduit la notion de \emph{syst\`eme de points} 
dans une dg-cat\'egorie satur\'ee $T$. Nous montrons que la donn\'ee
d'un tel syst\`eme permet de construire un espace alg\'ebrique $M_{\mathcal{P}}$, 
de type fini, lisse et s\'epar\'e, ainsi qu'un dg-foncteur de $T$ vers
une version tordue de la dg-cat\'egorie d\'eriv\'ee de $M_{\mathcal{P}}$.
On montre de plus que ce dg-foncteur est une \'equivalence si et seulement 
si $M_{\mathcal{P}}$ est propre. Tout au long de ce travail nous 
\'etudions les t-structures sur les familles alg\'ebriques d'objets 
dans $T$, ce qui poss\`ede possiblement un int\'er\^et en soi ind\'ependant
du th\`eme de ce travail.
\end{abstract}

\tableofcontents

\section*{Introduction}

Toute vari\'et\'e propre et lisse $X$ sur un corps $k$ (ici alg\'ebriquement clos) donne lieu 
\`a une dg-cat\'egorie $L_{parf}(X)$ des complexes parfaits sur $X$. 
La dg-cat\'egorie, souvent consid\'er\'ee comme l'espace non-commutatif
associ\'e \`a $X$, se rappelle de tr\`es nombreux invariants cohomologiques 
et g\'eom\'etriques de la vari\'et\'e $X$ (voir par exemple \cite{rouq,to2}). Cependant,
la pr\'esence de vari\'et\'es qui partagent un m\^eme cat\'egorie
d\'eriv\'ee implique qu'il est en g\'en\'eral impossible de reconstruire
$X$ \`a partir de $L_{parf}(X)$, et quand bien m\^eme une  
telle reconstruction est possible il existe de tr\`es nombreux choix pour $X$. 

Dans \cite{tova1} nous avons introduit une construction dans le sens inverse 
$T \mapsto \mathcal{M}_{T}$, qui \`a une dg-cat\'egorie $T$ associe 
le ($\infty$-)champ classifiant des objets dans $T$. Cette construction n'est pas
inverse de $X \mapsto L_{parf}(X)$, mais est adjointe en un sens
pr\'ecis (voir \cite[\S 3.1]{to4}). Lorsque $T$ s'\'ecrit de la forme $L_{parf}(X)$, 
la vari\'et\'e $X$ se retrouve comme un ouvert $X \subset \mathcal{M}_{T}$ (modulo
une $\mathbb{G}_{m}$-gerbe triviale), qui correspond \`a la partie de
$\mathcal{M}_{T}$ qui classe les faisceaux gratte-ciel associ\'es aux points de $X$
vus comme complexes parfaits sur $X$. Un autre choix
de vari\'et\'e $X'$ telle que $T \simeq L_{parf}(X')$ donne lieu \`a un
autre ouvert $X' \subset \mathcal{M}_{T}$. Ceci montre que 
la reconstruction d'une vari\'et\'e $X$ telle que $T\simeq L_{parf}(X)$
n'est envisageable que si l'on sp\'ecifie un sous-ensemble $\mathcal{P}$
des classes d'\'equivalence d'objets de $T$ correspondant
aux gratte-ciel de $X$ vus comme objets de $T$. 

Le r\'esultat principal de ce travail est le th\'eor\`eme
\ref{t2}, dans le quel nous donnons des conditions n\'ecessaires et suffisantes
sur un ensemble $\mathcal{P}$
de classes d'\'equivalence d'objets dans une dg-cat\'egorie $T$ propre et lisse, 
pour qu'il existe un espace alg\'ebrique $X$, lisse, s\'epar\'e et de type fini sur $k$, 
une classe 
$\alpha \in H^{2}_{et}(X,\mathbb{G}_{m})$
et un dg-foncteur 
$$\phi_{\mathcal{P}} : T^{op} \longrightarrow L_{parf}^{\alpha}(X),$$
qui identifie l'ensemble $\mathcal{P}$ aux gratte-ciel de $X$ (o\`u 
$L_{parf}^{\alpha}(X)$ est la dg-cat\'egorie des complexes parfaits sur
$X$ tordus par $\alpha$). On montre de plus que 
$\phi_{\mathcal{P}}$ est une \'equivalence si et seulement si 
$X$ est propre. Cet \'enonc\'e est une r\'eponse possible au probl\`eme
de savoir si une dg-cat\'egorie propre et lisse $T$ est d'origine 
g\'eom\'etrique. Cependant, plus qu'un simple th\'eor\`eme
de reconstruction nous pensons que notre r\'esultat 
peut \^etre un outil utile pour construire des \'equivalences
entre cat\'egories d\'eriv\'ees, bien que cet apsect ne sera pas discut\'e
en d\'etail dans ce travail (voir cependant \S 8). 

Les conditions sur la dg-cat\'egorie $T$ et l'ensemble
$\mathcal{P}$ qui constituent les hypoth\`eses de notre th\'eor\`eme
ne peuvent pas se r\'esumer dans cette introduction, et une grande partie
de ce travail consiste en l'introduction des notions qui entrent en jeu. 
Elles peuvent cependant se d\'ecliner en quatre grandes familles de conditions
que nous allons bri\`evement commenter.

\begin{enumerate}

\item La dg-cat\'egorie $T$ est satur\'ee.

\item Les objets de $\mathcal{P}$ sont des objets ponctuels deux \`a deux orthogonaux
et co-engendrent $T$.

\item Les objets de $\mathcal{P}$ co-engendrent une t-structure parfaite et ouverte 
sur $T$.

\item La famille des objets de $\mathcal{P}$ est born\'ee. 

\end{enumerate}

La premi\`ere condition affirme que $T$ est propre, lisse et triangul\'ee (voir \cite{tova1}), 
ce qui est une hypoth\`ese naturelle et incontournable si l'on souhaite
reconstruire une vari\'et\'e propre et lisse $X$. Par ailleurs, 
ces hypoth\`eses impliquent l'existence d'un foncteur de Serre pour $T$ (voir 
notre \S 1), ce qui intervient
de mani\`ere cruciale tout au long de l'article.

La condition $(2)$ contraint le comportement cohomologique
des objets de $\mathcal{P}$. On demande par exemple que l'on ait 
$T(x,x)\simeq Sym_{k}(Ext^1(x,x)[1])$ avec $Ext^1(x,x)$ de dimension uniforme $d$. 
Par ailleurs, si $x \neq y$ on demande que $T(x,y)\simeq 0$. On doit aussi
avoir $S_{T}(x)\simeq x[d]$ pour tout $x\in \mathcal{P}$, o\`u $S_T$ est
le foncteur de Serre de $T$.
Enfin, on demande que l'ensemble des objets de $\mathcal{P}$
soit une \emph{spanning class} au sens de Bridgeland:
si $T(y,x)\simeq 0$ pour tout $x \in \mathcal{P}$ alors $y\simeq 0$. 

La condition $(3)$ est la plus d\'elicate et la plus indirecte. Tout d'abord, 
on d\'eclare qu'un objet $y \in T$ est \emph{positif} si pour tout
$x\in \mathcal{P}$, le complexe $T(y,x)$ est cohomologiquement concentr\'e 
en degr\'e n\'egatif. Une premi\`ere condition est que cette notion
de positivit\'e d\'efinisse une t-structure sur la cat\'egorie
triangul\'ee $[T]$ associ\'ee \`a $T$. Cette t-structure est dite co-engendr\'ee
par l'ensemble $\mathcal{P}$. On demande de plus qu'elle soit \emph{ouverte}, 
c'est \`a dire qu'\^etre un objet du c\oe ur de cette structure soit une condition
ouverte (qui est une mani\`ere d'imposer une semi-continuit\'e
pour les objets de cohomologie associ\'es \`a la t-structure).
 Cette notion demande une \'etude des t-structures induites sur les familles
alg\'ebriques d'objets de $T$ (voir notre \S 3) et nous ne sommes pas arriv\'es \`a 
la d\'ecrire de mani\`ere simple en termes de $T$ seule, on encore
\`a donner des conditions suffisantes pour qu'elle soit automatique. On demande aussi 
que la t-structure soit \emph{parfaite}, ce qui signifie essentiellement que son
c\oe ur est noeth\'erien et que les foncteurs de troncations 
pr\'eservent les objets compacts. Encore une fois nous ne sommes pas arriv\'es
\`a d\'ecrire cette condition simplement en termes de $T$.

Enfin, la derni\`ere condition $(4)$ affirme que la famille d'objets
$\mathcal{P}$ vit dans une partie quasi-compacte du champ 
des objets $\mathcal{M}_T$ ce qui peut se traduire par 
l'existence d'une borne uniforme sur la cohomologie des objets 
de $\mathcal{P}$ par rapport \`a un g\'en\'erateur compact (voir \cite[\S 3.3]{tova1}).

Sous ces conditions $(1)$ \`a $(4)$, nous montrons l'existence d'un espace
de modules grossier $M_{\mathcal{P}}$ qui classe les objets de $\mathcal{P}$. 
On montre que cet espace est un espace alg\'ebrique de type fini, s\'epar\'e
et lisse sur $k$. Le champ des objets de $\mathcal{P}$ est lui une
$\mathbb{G}_{m}$-gerbe $\mathcal{M}_{\mathcal{P}} \longrightarrow M_{\mathcal{P}}$
d\'ecrit par une dg-alg\`ebre d'Azumaya (au sens de \cite{to4}) naturelle. 
Si on note $\alpha \in H^2_{et}(M_{\mathcal{P}},\mathbb{G}_{m})$ la classe
de cette gerbe, on montre l'existence d'un \emph{dg-foncteur de d\'ecomposition}
$$\phi_{\mathcal{P}} : T^{op} \longrightarrow L_{parf}^{\alpha}(M_{\mathcal{P}}),$$
qui est \`a rapprocher de la d\'ecomposition spectrale 
d'une repr\'esentation comme faisceau sur l'espace des repr\'esentations
irr\'eductibles d'un groupe. Ce dg-foncteur est une quasi-\'equivalence sur
les sous-dg-cat\'egories pleine form\'ees des points, et on montre que 
c'est une quasi-\'equivalence si et seulement si $M_{\mathcal{P}}$
est de plus propre (sans \^etre arriv\'e \`a mettre au jour un crit\`ere
simple qui implique la propret\'e de $M_{\mathcal{P}}$). \\

Pour conclure, nous souhaitons signaler que nous ne consid\'erons pas 
notre th\'eor\`eme principal comme r\'eellement optimal, mais plut\^ot comme
un premier pas, et qu'il 
serait souhaitable de l'am\'eliorer. Tout d'abord, comme nous l'avons
d\'ej\`a signal\'e, certaines conditions, particuli\`erement autour
de la t-structure, sont d\'elicates \`a d\'ecrire en termes
de $T$ seule. Leurs v\'erifications sont parfois difficiles dans 
la pratique, et il est fort probable que l'on puisse donner
des am\'eliorations (par exemple en donnant des conditions
suffisantes sur une t-structure pour \^etre ouverte ou encore parfaite). 
Par ailleurs, nous d\'emontrons que le foncteur 
$\phi_{\mathcal{P}} : T^{op} \longrightarrow L_{parf}^{\alpha}(M_{\mathcal{P}})$
est une \'equivalence seulement si $X$ est propre, sans donner pour autant 
de condition qui permette de s'en assurer. Encore une fois, dans
la pratique cela veut dire qu'il faut v\'erifier la propret\'e
de $X$ au cas par cas. On devrait pouvoir am\'eliorer cela 
en rajoutant, par exemple, des conditions de stabilit\'es en plus
de la donn\'ee de $\mathcal{P}$. A notre d\'echarge, notre point de d\'epart 
est une dg-cat\'egorie propre et lisse g\'en\'erale, qui peut ne rien avoir \`a voir, 
\`a priori, avec la dg-cat\'egorie d\'eriv\'ee d'une vari\'et\'e alg\'ebrique, 
on pourrait prendre par exemple des dg-cat\'egories d'origines
topologiques ou symplectiques (comme les cat\'egories de type Fukaya).
La conclucion du th\'eor\`eme \ref{t2} est ainsi relativement forte, et il n'est 
pas suprenant que les hypoth\`eses le soient proportionnellement.

Enfin, pour finir, 
nous avons inclus deux exemples de paire $(T,\mathcal{P})$ dans
notre dernier paragraphe, mais nous n'avons pas fait l'exercice
de v\'erifier toutes les conditions qui forment les hypoth\`eses
de notre th\'eor\`eme, notre but \'etant ici plus d'illustrer  
la signification du th\'eor\`eme plut\^ot que de donner de r\'eelles
applications. Il serait donc int\'eressant 
d'avoir un exemple d'application du th\'eor\`eme \ref{t2}
o\`u l'on obtienne une \'equivalence de cat\'egories d\'eriv\'ees qui n'\'etait
pas encore connue. Nous pensons par exemple \`a la sym\'etrie mirroir: 
une vari\'et\'e symplectique $N$ munie d'une fibration en tores lagrangiens $N\rightarrow S$,
semble pouvoir donner une paire $(T,\mathcal{P})$, o\`u $T$ est un dg-mod\`ele
pour la $\mathcal{A}_{\infty}$-cat\'egorie de Fukaya de $N$, et 
$\mathcal{P}$ est la famille d'objets correspondant \`a la famille de tores lagrangiens
donn\'ee par les fibres de $N \rightarrow S$. Ici, l'espace 
alg\'ebrique $M_{\mathcal{P}}$ serait bien entendu un candidat au miroir de $M$
(voir l'introduction de \cite{abou}), modulo
d'innombrables complications techniques (e.g. le corps de base $k$ doit \^etre remplac\'e
par quelque chose comme $k((t))$, etc).
Malheureusement, notre manque de compr\'ehension des cat\'egories de Fukaya
ne nous permet pas d'en dire plus dans ce travail. \\

\textbf{Conventions:} Tout au long de ce travail nous travaillerons 
au-dessus d'un corps alg\'ebriquement cl\^os $k$ de caract\'eristique nulle.

\section{Dg-cat\'egories satur\'ees}

Nous nous pla\c{c}ons dans le cadre de la th\'eorie Morita des dg-cat\'egories ($k$-lin\'eaires)
de \cite{to2}, auquel nous renvoyons le lecteur pour les d\'etails. Nous travaillerons
dans $Ho(dg-cat)$, la cat\'egorie homotopique des dg-cat\'egories $k$-lin\'eaires. Concr\`etement, 
cela signifie que nous travaillons \`a quasi-\'equivalence de dg-cat\'egories pr\`es, 
sans que cela soit dit de mani\`ere explicite. 
Il nous arrivera par exemple d'utiliser l'expression \emph{dg-foncteur} pour signifier
en r\'ealit\'e un morphisme dans $Ho(dg-cat)$, ou encore \emph{limites} pour 
signifier \emph{limites homotopiques}.  De m\^eme 
certaines constructions seront d\'ecrites de mani\`ere na\"\i ve
et demanderaient d'expliciter certains mod\`eles explicites (ce que nous
laissons au lecteur le soin de faire, ou pas). 
Nous rencontrerons parfois des dg-cat\'egories non petites, et nous laissons
le soins au lecteur de fixer des univers pour donner un sens 
\`a certains \'enonc\'es (voir par exemple \cite{to2,to4}). \\

Soit $T$ une petite dg-cat\'egorie sur $k$. Nous noterons 
$\widehat{T}:=\mathbb{R}\underline{Hom}(T^{op},\widehat{k})$, la dg-cat\'egorie
obtenue \`a partir de $T$ en y rajoutant des colimites homotopiques (voir \cite{to2}). 
L'objet $\widehat{T}$ est bien d\'efini dans $Ho(dg-cat)$, la cat\'egorie
homotopique des dg-cat\'egories. Un mod\`ele explicite de 
$\widehat{T}$ est la dg-cat\'egorie des $T^{op}$-dg-modules cofibrants et fibrants
(voir \cite[\S 2.2]{tova1}). 

On dispose d'un plongement de Yoneda $\underline{h} : T \longrightarrow 
\widehat{T}$ qui identifie $T$ \`a une sous-dg-cat\'egorie pleine
de $\widehat{T}$. Le morphisme $\underline{h}$ se factorise par $\widehat{T}_{c} \subset \widehat{T}$
la sous-dg-cat\'egorie pleine form\'ee des objets compacts. \\

Nous rappelons la terminologie suivante.

\begin{itemize}

\item La dg-cat\'egorie $T$ est \emph{triangul\'ee} si 
le morphisme $h$ induit une quasi-\'equivalence $T \simeq \widehat{T}_{c}$. 

\item La dg-cat\'egorie $T$ est \emph{lisse sur $k$} s'il existe une dg-alg\`ebre
$B$ sur $k$, avec $\widehat{T}$ quasi-\'equivalente \`a $\widehat{B}$, et telle
que $B$ soit un $B\otimes_{k}B^{op}$-dg-module compact. 

\item La dg-cat\'egorie $T$ est \emph{propre sur $k$} s'il existe une dg-alg\`ebre
$B$ sur $k$, avec $\widehat{T}$ quasi-\'equivalente \`a $\widehat{B}$, et telle
$B$ soit compacte comme $k$-dg-module (i.e. soit un complexe parfait 
de $k$-espaces vectoriels).

\item La dg-cat\'egorie $T$ est \emph{satur\'ee} si elle est 
propre, lisse et triangul\'ee. 

\end{itemize}

Les dg-cat\'egories satur\'ees poss\`edent de fabuleuses 
propri\'et\'es de dualit\'e. Elles sont par exemple les objets pleinement 
dualisables de la cat\'egorie mono\"\i dale sym\'etrique des petites dg-cat\'egories
\`a \'equivalence de Morita pr\`es (voir \cite[Prop. 2.5]{to4}), et aussi les
objets pleinement dualisables de la $2$-cat\'egorie mono\"\i dale sym\'etrique des
dg-cat\'egories compactement engendr\'ees et dg-foncteurs continus (voir \cite{lu2}). 
En ce qui nous concerne, nous rentiendrons les faits suivants.

\begin{enumerate}

\item Pour deux dg-cat\'egories satur\'ees $T_1$ et $T_2$, 
la dg-cat\'egorie $\mathbb{R}\underline{Hom}(T_1,T_2)$ est 
satur\'ee. Elle est de plus naturellement \'equivalente
\`a $\widehat{(T_1\otimes T_2^{op})}_{c}$, la dg-cat\'egorie
des $(T_1,T_2)$-dg-bimodules compacts.  Le dg-bimodule correspondant 
\`a un morphisme $f$ sera g\'en\'eralement not\'e $\Gamma(f)$. 
Ce dg-bimodule envoie $(a,b) \in T_1 \otimes T_2^{op}$ sur le 
complexe $T_2(b,f(a))$. 

\item Toute dg-cat\'egorie satur\'ee $T$ poss\`ede un 
foncteur de Serre $S_{T}$. Il s'agit d'un endomorphisme 
$S_{T} : T \longrightarrow T$ de $T$, tel que le dg-bimodule correspondant soit 
donn\'e par la formule
$$\Gamma(S_{T}) : (a,b) \mapsto T(a,b)^{\vee}:=\underline{Hom}(T(a,b),k).$$
On dispose donc de quasi-isomorphisme naturel
$$T(a,b)^{\vee} \simeq T(b,S_{T}(a)).$$
L'endomorphisme $S_T$ de $T$ est toujours une auto-\'equivalence. 

\item Soit $f : T_1 \longrightarrow T_2$ un dg-foncteur entre 
dg-cat\'egories satur\'ees. Alors $f$ poss\`ede un adjoint \`a droite
et un adjoint \`a gauche (au sens de \cite[\S 3.1]{to4}), not\'es respectivement
$$f_* : T_2 \longrightarrow T_1 \qquad  f_! : T_2 \longrightarrow T_1.$$
Si $f$ correspond au dg-bimodule $\Gamma(f) \in \widehat{T_1\otimes T_2^{op}}$, alors
$f_*$ correspond au $(T_2,T_1)$-dg-bimodule
$$\Gamma(f_!) : (a,b) \mapsto T_2(f(b),a).$$ 
L'adjoint \`a gauche $f_!$ est quand \`a lui donn\'e par la formule
$$f_! = S_{T_1}^{-1}(f_*)S_{T_2},$$
o\`u $S_{T_{i}}$ est le foncteur de Serre de la dg-cat\'egorie $T_i$.
En d'autres termes le dg-bimodule correspondant \`a $f_!$
est $(a,b) \mapsto T_2(a,fS_{T_1}(b))^{\vee}$.

\end{enumerate}

\section{Espaces de modules d'objets simples dans une dg-cat\'egorie satur\'ee}

Soit $T$ une dg-cat\'egorie satur\'ee sur $k$. On d\'efinit un pr\'efaisceau simplicial
sur la cat\'egorie $Aff_k$ des sch\'emas affines sur $k$
$$\mathcal{M}_{T} : Aff_k^{op} \longrightarrow sSet,$$
par la formule
$$\mathcal{M}_{T}(Spec\, A):=Map_{dg-cat}(T^{op},\widehat{A}_{c}).$$
Le pr\'efaisceau $\mathcal{M}_{T}$ est un champ pour la topologie fpqc (voir \cite{tova1}). 
De plus, pour tout champ $F \in \St_{k}$ on dispose d'une \'equivalence naturelle
$$Map_{\St_{k}}(F,\mathcal{M}_{T}) \simeq Map_{dg-cat}(T^{op},L_{parf}(F)),$$
o\`u $L_{parf}(F)$ est la dg-cat\'egorie des complexes parfaits sur 
$F$ d\'efinie par
$$L_{parf}(F):=\lim_{Spec\, A \longrightarrow F} \widehat{A}_{parf}.$$
Il s'agit ici d'une limite homotopique dans la th\'eorie homotopique
des dg-cat\'egories (voir \cite{hagII} pour plus de d\'etails sur la descente).

Le principal th\'eor\`eme de \cite{tova1} affirme que $\mathcal{M}_{T}$ est un 
champ localement alg\'ebrique et localement de pr\'esentation finie sur $k$. 
Dans ce travail nous nous restreindrons au sous-champ $\mathcal{M}_{T}^{simp} \subset 
\mathcal{M}_{T}$, form\'e des objets simples. Par d\'efinition, un 
point de $\mathcal{M}_{T}(Spec\, A)$ correspond \`a un $T^{op}\otimes_{k}A$-dg-module
compact. Un tel $T^{op}\otimes_{k}A$-dg-module
$E$ poss\`ede un complexe d'endomorphismes $\mathbb{R}\underline{Hom}(E,E)$ 
qui est un complexe de $A$-dg-modules. Ce complexe est de plus parfait sur $A$ car 
$T$ est satur\'ee. Par d\'efinition, l'objet $E \in \mathcal{M}_{T}(Spec\, A)$
est dit \emph{simple} si le complexe parfait $\mathbb{R}\underline{Hom}(E,E)$ 
v\'erifie les deux conditions suivantes.

\begin{itemize}

\item L'amplitude de $\mathbb{R}\underline{Hom}(E,E)$  est 
contenue dans $[0,\infty]$. 

\item Le morphisme naturel $A \longrightarrow \mathbb{R}\underline{Hom}(E,E)$ 
induit par l'identit\'e de $E$, induit pour tout morphisme d'anneaux 
$A \longrightarrow A'$ un isomorphisme
$$A' \simeq H^{0}(\mathbb{R}\underline{Hom}(E,E\otimes_{A}A')).$$
\end{itemize}

Par semi-continuit\'e on voit que l'inclusion naturelle $\mathcal{M}_{T}^{simp} \subset 
\mathcal{M}_{T}$ est une immersion ouverte de champs. On en d\'eduit donc
que $\mathcal{M}_{T}^{simp}$ est lui-m\^eme un champ localement alg\'ebrique et localement
de pr\'esentation finie. Par ailleurs, l'annulation des $Ext$ n\'egatifs des objets simples
implique que $\mathcal{M}_{T}^{simp}$ est un $1$-champ d'Artin localement de pr\'esentation finie 
(voir \cite[\S 3.4]{tova1}). 
Enfin, le fait que les objets simples ne poss\`edent que les scalaires comme
endomorphismes de degr\'e z\'ero implique que $\mathcal{M}_{T}^{simp}$ est 
une $\mathbb{G}_{m}$-gerbe au-dessus d'un espace alg\'ebrique $M_{T}^{simp}$
localement de pr\'esentation finie sur $k$. De mani\`ere plus explicite, 
la projection
$$\pi : \mathcal{M}_{T}^{simp} \longrightarrow M_{T}^{simp}$$
est le morphisme quotient pour l'action naturelle du champ en groupes 
$K(\mathbb{G}_{m},1)$ sur $\mathcal{M}_{T}^{simp}$ qui consiste 
\`a tensoriser par des fibr\'es en droites de la base. Pour $Spec\, A \in Aff_k$, 
le groupe simplicial $K(A^{*},1)$ op\`ere naturellement sur 
l'ensemble simplicial $\mathcal{M}_{T}(A)$ par son action naturelle
sur la dg-cat\'egorie $\widehat{A}_{c}$ des complexes parfaits de $A$-dg-modules. 
Lorsque $Spec\, A$ d\'ecrit $Aff_k$ cela d\'efinit une action du champ 
en groupes $K(\mathbb{G}_{m},1)$ sur $\mathcal{M}^{simp}_{T}$ dont le quotient 
est $M_{T}^{simp}$. On peut 
donc \'ecrire
$$ M_{T}^{simp} \simeq [\mathcal{M}_{T}^{simp}/K(\mathbb{G}_{m},1)]$$
de sorte \`a ce que la projection $\pi : \mathcal{M}_{T}^{simp} \longrightarrow M_{T}^{simp}$
soit le morphisme quotient. 

Notons que la gerbe $\mathcal{M}_{T}^{simp}$ au-dessus
de $M_{T}^{simp}$ est en g\'en\'eral non-triviale. Elle est cependant donn\'ee par 
une dg-alg\`ebre d'Azumaya sur $M_{T}^{simp}$ au sens de \cite{to4}, qui peut-\^etre
construite de la mani\`ere suivante. Sur le champ 
$\mathcal{M}_{T}^{simp}$, on dispose d'une $T^{op}$-dg-module universel
$\mathcal{E}$, qui apr\`es oubli de la structure de $T^{op}$-dg-module 
fournit un complexe parfait $\mathcal{E}_{0}$ sur $\mathcal{M}_{T}^{simp}$. 
La dg-cat\'egorie des complexes parfaits sur $\mathcal{M}_{T}^{simp}$ 
poss\`ede des sous-dg-cat\'egories pleines
$$L_{parf}^{\chi}(\mathcal{M}_{T}^{simp}) \subset L_{parf}(\mathcal{M}_{T}^{simp}),$$
o\`u $\chi$ parcours l'ensemble des caract\`eres $\mathbb{Z}$ de $\mathbb{G}_{m}$, 
et o\`u $L_{parf}^{\chi}(\mathcal{M}_{T}^{simp})$ consiste en les complexes parfaits
de poids $\chi$: un complexe parfait $E \in L_{parf}(\mathcal{M}_{T}^{simp})$
vit dans $L_{parf}^{\chi}(\mathcal{M}_{T}^{simp})$ si et seulement si 
pour toute gerbe r\'esiduelle $B\mathbb{G}_{m} \subset \mathcal{M}_{T}^{simp}$
la restriction de $E$ \`a $B\mathbb{G}_{m}$ est un complexe de
repr\'esentations de $\mathbb{G}_{m}$ pures de poids $1$. Notons que
si $\mathcal{M}_{T}^{simp,\nu}$ d\'esigne un ouvert quasi-compact
de $\mathcal{M}_{T}^{simp}$ (voir \cite[\S 3.3]{tova1}), 
alors les inclusions naturelles d\'efinissent 
une quasi-\'equivalence de dg-cat\'egories
$$\bigoplus_{\chi \in \mathbb{Z}} L_{parf}^{\chi}(\mathcal{M}_{T}^{simp,\nu})
\simeq L_{parf}(\mathcal{M}_{T}^{simp,\nu}).$$
Cette d\'ecomposition n'est plus valable sur $\mathcal{M}_{T}^{simp}$
tout entier d\^ue \`a la non-quasi-compacit\'e. 

Nous disposons donc d'un complexe parfait universel $\mathcal{E}_{0}$ sur
$\mathcal{M}_{T}^{simp}$, pur de poids $1$. Cet objet est de plus un g\'en\'erateur
compact \emph{local} de $L_{parf}^{\chi=1}(\mathcal{M}_{T}^{simp})$, au sens
o\`u son image r\'eciproque sur tout affine $Spec\, A \longrightarrow \mathcal{M}_{T}^{simp}$
est un g\'en\'erateur compact de $D(A)$ la cat\'egorie d\'eriv\'ee de $A$. 
Il s'en suit, d'apr\`es le formalisme g\'en\'eral de \cite{to4}, que 
le faisceau en dg-alg\`ebres parfaites $\mathcal{A}:=\mathbb{R}\underline{End}(\mathcal{E}_{0})$
est pur de poids $0$, donc vit dans 
$L^{\chi=0}_{parf}(\mathcal{M}_{T}^{simp})\simeq L_{parf}(M_{T})$, et est de plus une dg-alg\`ebre
d'Azumaya sur l'espace alg\'ebrique
$M_{T}^{simp}$. D'apr\`es \cite[Cor. 4.8]{to4}, elle d\'etermine une classe $\gamma(\mathcal{A}) \in 
H^{2}_{et}(M_{T}^{simp},\mathbb{G}_{m})$ qui n'est autre que la classe
de la gerbe $\mathcal{M}_{T}^{simp} \longrightarrow M_{T}^{simp}$. Par construction, on a une
quasi-\'equivalence de dg-cat\'egories
$$L_{parf}(\mathcal{A}) \simeq L_{parf}^{\chi=1}(\mathcal{M}_{T}^{simp}),$$
o\`u le membre de gauche est la dg-cat\'egorie des $\mathcal{A}$-dg-modules parfaits sur
$M_{T}^{simp}$. Nous renvoyons \`a \cite{to4} pour plus de d\'etails.

\section{Extension des scalaires et t-structures}

Comme pr\'ec\'edemment, soit $T$ une dg-cat\'egorie triangul\'ee sur $k$ et 
$\widehat{T}$ la dg-cat\'egorie des dg-modules sur $T^{op}$. On note
$\underline{h} : T \longrightarrow \widehat{T}$ le dg-foncteur
de Yoneda, qui est un dg-foncteur pleinement fid\`ele. On note $[T]:=H^{0}(T)$
la cat\'egorie homotopique de $T$. Comme $T$ est triangul\'ee, $[T]$ est  
munie d'une structure triangul\'ee naturelle pour laquelle les triangles distingu\'es sont
les images des suites exactes de fibrations (voir \cite[\S 7]{ho} ou \cite[\S 4.4]{to3}). D'apr\`es 
nos conventions 
$[T]$ est aussi Karoubienne. De la m\^eme mani\`ere, $[\widehat{T}]$ est une
cat\'egorie triangul\'ee qui poss\`ede des sommes infinies. 
Le plongement de Yoneda induit un foncteur 
triangul\'e pleinement fid\`ele
$$[T] \hookrightarrow [\widehat{T}].$$
Le caract\`ere triangul\'e de $T$ implique que l'image essentielle de
ce plongement est la sous-cat\'egorie pleine des objets compacts 
dans $[\widehat{T}]$. On dipose donc d'une \'equivalence triangul\'ee naturelle
$$[T] \simeq [\widehat{T}]_{c}.$$

\begin{df}\label{d1}
Une \emph{t-structure sur $T$} est la donn\'ee d'une t-structure 
sur la cat\'egorie triangul\'ee $[\widehat{T}]$ au sens de \cite{bbd}. 
Le \emph{c\oe ur} d'une t-structure sur $T$ est par d\'efinition 
le c\oe ur de la t-structure correspondante sur $[\widehat{T}]$. Il sera
not\'e $\widehat{\mathcal{H}}$. 
\end{df}

Par d\'efinition, une t-structure sur $T$ est compl\`etement caract\'eris\'ee
par la donn\'ee d'une sous-cat\'egorie pleine $[\widehat{T}]^{\leq 0} \subset 
[\widehat{T}]$, form\'ee des objets $x$ tels que $\tau_{> 0}x\simeq 0$. 
Inversement, une telle sous-cat\'egorie pleine de $[\widehat{T}]$ d\'etermine
une t-structure si et seulement si elle est stable par sommes (possiblement infinies), 
c\^ones et facteurs directs, et si de plus elle est engendr\'ee
par c\^ones, sommes et facteurs directs par un ensemble petit d'objets
(voir \cite{kv}). Par la suite, nous supposerons implicitement que 
les t-structures consid\'er\'ees satisferont toutes cette hypoth\`eses
ensemblistes: $[\widehat{T}]^{\leq 0}$ est engendr\'ee, par 
c\^ones, sommes et facteurs directs par un ensemble petit d'objets.

Comme l'ensemble des classes 
d'\'equivalence de sous-cat\'egories pleines 
de $[\widehat{T}]$ sont en correspondance bi-univoque avec celui des sous-dg-cat\'egories
pleine de $\widehat{T}$, on voit qu'une t-structure sur $T$ consiste en la donn\'ee
d'une sous-dg-cat\'egorie pleine $\widehat{T}^{\leq 0} \subset \widehat{T}$
qui est stable, \`a \'equivalence pr\`es, par colimites (au sens des
dg-cat\'egorie) et qui est engendr\'ee par colimites par un ensemble petit d'objets.
Il est important de noter cependant que bien que l'inclusion
naturelle $[\widehat{T}]^{\leq 0} \subset [\widehat{T}]$ soit induite 
par un dg-foncteur pleinement fid\`ele $\widehat{T}^{\leq 0} \subset \widehat{T}$, 
le foncteur de troncation
$$\tau_{\leq 0} : [\widehat{T}] \longrightarrow [\widehat{T}]^{\leq 0}$$
ne l'est pas (car il n'est pas un foncteur triangul\'e). Ce foncteur
de troncation peut cependant \^etre repr\'esent\'e par un $\s$-foncteur
sur les $\s$-cat\'egories correspondantes
$$|\widehat{T}| \longrightarrow |\widehat{T}^{\leq 0}|$$
(o\`u $|-|$ d\'esigne le foncteur qui \`a une dg-cat\'egorie
associe l'$\s$-cat\'egorie correspondante, en appliquant par exemple
la construction de Dold-Kan aux complexes de morphismes, voir par exemple
\cite{ta3}). \\

Fixons maintenant une dg-cat\'egorie $T$ triangul\'ee sur $k$ munie d'une
t-structure. Nous noterons $\widehat{T}^{\leq 0}, \widehat{T}^{\geq 0} \subset \widehat{T}$
les sous-dg-cat\'egories pleines form\'ees des objets n\'egatifs, respectivement positifs,
par rapport \`a la t-structure. Pour $n\in \mathbb{Z}$, on pose
$\widehat{T}^{\leq n} \subset \widehat{T}$ la sous-dg-cat\'egorie pleine
form\'ee des objets $x$ tels que $\tau_{> 0}(x[n])\simeq 0$. De m\^eme, 
$\widehat{T}^{\geq n}$ est d\'efinie comme la sous-dg-cat\'egorie pleine form\'ee
des $x$ tels que $\tau_{<0}(x[n])\simeq 0$. 
On pose, pour toute paire d'entiers $a\leq b$
$$\widehat{T}^{[a,b]} := \widehat{T}^{\leq b} \cap \widehat{T}^{\geq a} \subset 
\widehat{T}.$$
Notons que $\widehat{T}^{[0,0]}$ est une sous-dg-cat\'egorie
pleine de $\widehat{T}$ dont la cat\'egorie homotopique 
s'identifie naturellement au c\oe ur de la t-structure. De plus, 
pour deux objets $x,y\in \widehat{T}^{[0,0]}$ les complexes
de morphismes $T(x,y)$ sont cohomologiquement concentr\'es en 
degr\'es positifs. On voit ainsi que $[\widehat{T}^{[0,0]}]$, qui est une cat\'egorie
ab\'elienne $k$-lin\'eaire, est munie d'une dg-foncteur naturel
$$[\widehat{T}^{[0,0]}] \longrightarrow \widehat{T}^{[0,0]}$$
induisant une \'equivalence sur les cat\'egorie homotopiques. En composant 
avec l'inclusion naturelle on trouve un dg-foncteur 
$$[\widehat{T}^{[0,0]}] \longrightarrow \widehat{T}.$$
Ce dg-foncteur n'est pas pleinement fid\`ele, mais induit un foncteur pleinement fid\`ele
$$[\widehat{T}^{[0,0]}] \hookrightarrow [\widehat{T}]$$
qui identifie $[\widehat{T}^{[0,0]}]$ au c\oe ur. De cette mani\`ere, le c\oe ur
de la t-structure $\widehat{\mathcal{H}}$ sera consid\'er\'e 
comme sous-cat\'egorie pleine de $[\widehat{T}]$, mais aussi comme
sous-dg-cat\'egorie (non-pleine) de $\widehat{T}$. \\

Nous introduisons les notions de finitudes suivantes pour une t-structure
sur $T$.

\begin{df}\label{d1'}
Soit $T$ une dg-cat\'egorie triangul\'ee sur $k$ munie d'une t-structure. 
\begin{enumerate}

\item Nous dirons que la t-structure est \emph{born\'ee} si 
tous les objets compacts de $[\widehat{T}]$ sont born\'es: pour tout
$E \in [T]\simeq [\widehat{T}]_{c}$ il existe $a\leq b$ tels que $E \in 
\widehat{T}^{[a,b]}$.

\item Nous dirons que la t-structure est \emph{pr\'ecompacte} si 
pour tout $E \in [T]$, les objets $\tau_{\leq 0}(E)$ et $\tau_{\geq 0}(E)$
sont compacts dans $[\widehat{T}]$.

\item Nous dirons que la t-structure est \emph{compacte} si elle
est born\'ee, pr\'ecompacte et si de plus la sous-dg-cat\'egorie
pleine $\widehat{T}^{\geq 0} \subset \widehat{T}$ est stable
par colimites filtrantes. 

\end{enumerate}
\end{df}

Les t-structures compactes d\'efinies ci-dessus poss\`edent des 
propri\'et\'es de stabilit\'e remarquables que nous avons
rassembl\'ees dans la proposition suivante.

\begin{prop}\label{p0}
Soit $T$ une dg-cat\'egorie triangul\'ee munie d'une t-structure
compacte. Les propri\'et\'es suivantes sont satisfaites.
\begin{enumerate}

\item Pour toute paire d'entiers $a\leq b$, la sous-dg-cat\'egorie
$\widehat{T}^{[a,b]} \subset \widehat{T}$ est stable par colimites
filtrantes. Les $\s$-foncteurs de cohomologie $H^{i}_{t} : |\widehat{T}| \longrightarrow
\widehat{\mathcal{H}}$ commutent aux colimites filtrantes.

\item Pour tout $-\infty\leq a\leq b \leq \infty$, la sous-dg-cat\'egorie
pleine $\widehat{T}^{[a,b]} \subset \widehat{T}$ est engendr\'ee par 
colimites par les objets de $T \cap \widehat{T}^{[a,b]}$. 

\item La t-structure sur $[\widehat{T}]$ se restreint en une $t$-structure
sur la sous-cat\'egorie des objets compacts $[T]$. Le c\oe ur $\mathcal{H}$ de cette 
t-structure induite est $\widehat{\mathcal{H}}\cap [T] \subset [\widehat{T}]$, l'intersection 
du c\oe ur de la t-structure sur $[\widehat{T}]$ avec la sous-cat\'egorie des
objets compacts.

\item La sous-cat\'egorie des objets $\omega$-petits de la cat\'egorie
ab\'elienne $\widehat{\mathcal{H}}$ est $\mathcal{H}=\widehat{\mathcal{H}}\cap [T]$. 

\item Un objet $E \in \widehat{T}$ est compact si et seulement 
s'il existe deux entiers $a\leq b$ avec $E \in \widehat{T}^{[a,b]}$ et
de plus $H^{i}_{t}(E) \in \widehat{\mathcal{H}}$ est un objet de $\mathcal{H}$.

\end{enumerate}
\end{prop}

\textit{Preuve:} $(1)$ Pour commencer, notons $|\widehat{T}|$ l'$\infty$-cat\'egorie
associ\'ee \`a la dg-cat\'egorie $\widehat{T}$. Par d\'efinition 
de t-structure compacte la sous-$\infty$-cat\'egorie pleine 
$|\widehat{T}^{\geq 0}| \hookrightarrow |\widehat{T}|$ est stable by colimites filtrantes.
Comme l'$\infty$-foncteur de troncation $\tau_{\geq 0}$ est un adjoint \`a gauche, il 
commute aux colimites filtrantes. En utilisant le triangle distingu\'e 
d'$\infty$-foncteurs
$$\xymatrix{ \tau_{\leq - 1} \ar[r] & id \ar[r] & \tau_{\geq 0}}$$
on voit que $\tau_{\leq -1}$ commute aussi aux colimites filtrantes. Cela implique
que pour tout $a\leq b$ les $\s$-foncteurs $\tau_{[a,b]}$ commutent aux colimites fitrantes, et donc
que les sous $\s$-cat\'egories pleines $|\widehat{T}^{[a,b]}| \hookrightarrow 
|\widehat{T}|$ sont stables par colimites filtrantes. 

$(2)$ Soit $x\in \widehat{T}^{[a,b]}$. On \'ecrit $x\simeq \colim_{\alpha} x_{\alpha}$, 
une colimite filtrante d'objets $x_\alpha \in T$. Par $(1)$ on a 
$x \simeq \tau_{[a,b]}(x) \simeq \colim_{\alpha}\tau_{[a,b]}(x_{\alpha})$. 
Or, par hypoth\`ese de compacit\'e de la t-structure, 
tous les $\tau_{[a,b]}(x_{\alpha})$ sont dans $T$, et donc dans $T\cap \widehat{T}^{[a,b]}$. 

$(3)$ C'est une cons\'equence directe de la condition de compacit\'e sur 
la t-structure. 

$(4)$ Tout d'abord, comme les objets de $\mathcal{H}$ sont compacts dans $\widehat{T}$, 
ils sont aussi $\omega$-petits dans $\widehat{\mathcal{H}}$ (car $(1)$ implique que les colimites
filtrantes dans $\widehat{\mathcal{H}}$ sont aussi des colimites filtrantes dans $\widehat{T}$). 
Inversement, soit $x$ un objet $\omega$-petit de $\widehat{\mathcal{H}}$. On voit $x$ dans 
$\widehat{T}$ et on \'ecrit $x$ comme une colimite filtrante 
$x\simeq \colim_{\alpha} x_{\alpha}$ avec 
$x_{\alpha} \in T\cap \widehat{\mathcal{H}} = \mathcal{H}$ (ce qui est possible
\`a l'aide de $(2)$). Tous les objets $x$ et $x_{\alpha}$ \'etant $\omega$-petits dans 
$\widehat{\mathcal{H}}$ on voit que $x$ est un r\'etracte d'un $x_{\alpha}$. 
Mais $T\cap \mathcal{H}=\widehat{\mathcal{H}}$ est stable par facteurs directs dans $\widehat{T}$, 
et donc $x \in \mathcal{H}$. 

$(5)$ Se d\'eduit par d\'evissage des points pr\'ec\'edents.
\hfill $\Box$ \\

Soit maintenant $A$ une $k$-alg\`ebre commutative et consid\'erons
$T_{A}:=T\hat{\otimes}_{k}A$ la dg-cat\'egorie triangul\'ee $A$-lin\'eaire induite
par changement de base. La construction $A \mapsto T_{A}$ se promeut en 
un $\s$-foncteur de la cat\'egorie des $k$-alg\`ebres
commutatives vers celle des dg-cat\'egories triangul\'ees (voir \cite[\S 3.1]{to4}). 
De mani\`ere explicite, $T_{A}$ peut-\^etre d\'efinie, \`a \'equivalence naturelle pr\`es, 
comme la sous-dg-cat\'egorie pleine des $(T\otimes_{k}A)^{op}$-dg-modules, form\'ee
des dg-modules cofibrants et quasi-repr\'esentables (voir \cite{tova1}). On peut aussi 
voir $T_A$ par la formule suivante
$$T_{A}\simeq \mathbb{R}\underline{Hom}(A,T),$$
o\`u $A$ est consid\'er\'ee comme une dg-cat\'egorie sur $k$ avec un unique objet. 
La dg-cat\'egorie $\widehat{T_{A}}$ s'identifie naturellement 
\`a $\widehat{T}_{A}=\mathbb{R}\underline{Hom}(A,\widehat{T})$ ainsi qu'\`a la
dg-cat\'egorie des $(T\otimes_{k}A)^{op}$-dg-modules cofibrants.
Nous renvoyons \`a \cite{to4} pour plus de d\'etails sur les changements d'anneaux de base
pour les dg-cat\'egories. 

On d\'efinit une t-structure sur $\widehat{T}_{A}$ induite par celle de $T$ de la mani\`ere
suivante. On dispose d'un unique dg-foncteur  $k \longrightarrow A$
induit par l'identit\'e dans $A$. Cela d\'efinit un dg-foncteur de restriction
$\mathbb{R}\underline{Hom}(A,\widehat{T}) \longrightarrow 
\mathbb{R}\underline{Hom}(k,\widehat{T})\simeq \widehat{T}$, 
et donc un dg-foncteur, que nous appellerons \emph{oubli}
$$\widehat{T}_A \longrightarrow \widehat{T}.$$
On d\'efinit $\widehat{T}_{A}^{\geq 0}$  par la carr\'e cart\'esien suivant
$$\xymatrix{
\widehat{T}_{A}^{\geq 0} \ar[r] \ar[d] & \widehat{T}_{A} \ar[d] \\
\widehat{T}^{\geq 0} \ar[r] & \widehat{T}.}$$
Le dg-foncteur d'oubli admet un adjoint \`a gauche (au 
sens des dg-cat\'egories, voir \cite{to4})
$$-\otimes_{k}A : \widehat{T} \longrightarrow \widehat{T}_{A},$$
qui en termes de dg-modules envoie un $T^{op}$-dg-module $E$ sur
$E\otimes_{k}A$ muni de sa structure naturelle de $T^{op}\otimes_{k}A$-dg-module. 
On en d\'eduit un adjoint \`a gauche au foncteur d'oubli sur les cat\'egories
triangul\'ees associ\'ees, que nous appellerons \emph{changement de base}
$$-\otimes_{k}A : [\widehat{T}] \longrightarrow [\widehat{T}_A].$$
On d\'efinit alors $\widehat{T}_{A}^{\leq 0}$ comme \'etant la plus petite sous-dg-cat\'egorie
pleine de $\widehat{T}_{A}$ qui est stable par sommes, c\^ones et r\'etractes, et qui 
contient l'image de $\widehat{T}^{\leq 0}$ par le dg-foncteur $-\otimes_{k}A$.

La paire de sous-cat\'egories de $[\widehat{T}_{A}]$
$$[\widehat{T}_{A}]^{\leq 0}:=[\widehat{T}_{A}^{\leq 0}] \qquad 
[\widehat{T}_{A}]^{\geq 0}:=[\widehat{T}_{A}^{\geq 0}]$$
d\'efinit une t-structure sur $\widehat{T}_{A}$. Pour cette
t-structure, l'adjonction de foncteurs triangul\'es
$$-\otimes_{k}A : [\widehat{T}] \longleftrightarrow [\widehat{T}_A]$$
est compatible avec les t-structures au sens o\`u $-\otimes_{k}A$
est t-exact \`a droite, et son adjoint \`a droite t-exact \`a gauche. 

\begin{df}\label{d2}
Avec les notations ci-dessus, la t-structure sur $\widehat{T}_A$ est 
appel\'ee la \emph{t-structures induite par extension des scalaires}. 
\end{df}

Le lemme suivant se d\'eduit ais\'ement des d\'efinitions et du fait que
le foncteur d'oubli $[\widehat{T}_{A}] \longrightarrow [\widehat{T}]$ 
commute avec les sommes arbitraires.

\begin{lem}\label{l1}
Avec les notations ci-dessus, pour toute dg-cat\'egorie triangul\'ee
$T$ munie d'une t-structure, et toute $k$-alg\`ebre
commutative $A$, le foncteur d'oubli
$$p : [\widehat{T}_A] \longrightarrow [\widehat{T}]$$
est t-exact pour la t-structure sur $\widehat{T}_A$ induite par extension
des scalaires. De plus, ce foncteur refl\^ete la t-structure sur
$[\widehat{T}_A]$: pour tout objet $x \in [\widehat{T}_A]$ on a
$$( p(x) \in [\widehat{T}]^{\leq 0} ) \Longleftrightarrow
( x \in [\widehat{T}_A]^{\leq 0} ) \qquad ( p(x) \in [\widehat{T}]^{\geq 0}) \Longleftrightarrow
( x \in [\widehat{T}_A]^{\geq 0}).$$
\end{lem}

\textit{Preuve:} Il suffit de d\'emontrer les deux \'equivalences. 
La seconde \'equivalence est vraie par d\'efinition de $[\widehat{T}_A]^{\geq 0}$
comme image r\'eciproque de $[\widehat{T}]^{\geq 0}$ par $p$. Soit $x \in \widehat{T}_A^{\leq 0}$. 
C'est une colimite dans $\widehat{T}_A$ d'objets de la forme $x_\alpha \otimes A$, avec 
$x_\alpha \in \widehat{T}^{\leq 
0}$. Comme le dg-foncteur $p$ commute aux colimites arbitraires, $p(x)$ est 
une colimite d'objets de la forme $x_\alpha \otimes_{k} A$, qui eux m\^eme
sont des sommes de $x_\alpha$. On voit donc que $p(x) \in \widehat{T}^{\leq 0}$. Inversement, 
supposons que $p(x) \in \widehat{T}^{\leq 0}$. A l'aide de l'adjonction $(\otimes_k A,p)$, 
on construit un objet simplicial $R_*(x)$ dans $\widehat{T}_A$ muni
d'une augmentation $R_*(x) \longrightarrow x$ avec
$R_n(x)\simeq q^{n}(x)\simeq p(x)\otimes_{k} A^{\otimes n+1}$ pour tout $n\geq 0$, 
avec $q$ l'endomorphisme de $\widehat{T}_A$ donn\'e par $p(-)\otimes_{k}A$.  
Cet objet est une r\'esolution de $x$ au sens o\`u l'augmentation induit
une \'equivalence dans $\widehat{T}_A$
$$\colim_{n\in \Delta^{op}}R_n(x) \simeq x$$
Par hypoth\`ese et par d\'efinition chaque $R_n(x)$ est dans $\widehat{T}_A^{\leq 0}$,
et ainsi $x\in \widehat{T}_A^{\leq 0}$. 
\hfill $\Box$ \\

Ce lemme, associ\'e au fait que 
$\widehat{T}_A$ est la dg-cat\'egorie des $A$-modules dans $\widehat{T}$, montre que
le c\oe ur de la t-structure induite sur $\widehat{T}_A$ n'est autre que 
la cat\'egorie des $A$-module dans la cat\'egorie ab\'elienne $k$-lin\'eaire
$\widehat{\mathcal{H}}$. En d'autre termes, le c\oe ur de 
$\widehat{T}_A$ s'identifie \`a la cat\'egorie ab\'elienne $A$-lin\'eaire
obtenue \`a partir de $\widehat{\mathcal{H}}$ par extension des scalaires 
de $k$ \`a $A$. \\

Nous arrivons aux notions principales de cette section, la notion 
de \emph{t-structure ouverte}  et de \emph{t-structure parfaite} 
sur une dg-cat\'egorie satur\'ee $T$. Ces notions font intervenir 
les t-structures induites sur $\widehat{T}_A$ et ne semblent pas
facilement exprimables en terme de $T$ seule. 

Commen\c{c}ons par la notion de t-structure ouverte. 
Pour toute paire d'entiers $a \leq b$ nous 
d\'efinissons $\mathcal{M}_{T}^{[a,b]} \subset \mathcal{M}_{T}$
un sous-pr\'echamp de $\mathcal{M}_{T}$ comme suit. Soit $A$ une
$k$-alg\`ebre commutative, rappelons que l'ensemble simplicial
$\mathcal{M}_{T}(A)$ est le nerf de la cat\'egorie
des $T^{op}\otimes_{k}A$-dg-modules cofibrants et compact et des quasi-isomorphismes
entre iceux (voir \cite{tova1}). Par d\'efinition, $\mathcal{M}^{[a,b]}_{T}(A) \subset 
\mathcal{M}_{T}(A)$ est le sous-ensemble simplicial plein (i.e. 
r\'eunion de composantes connexes) form\'e des dg-modules $E$
satisfaisant la condition suivante: pour tout $A$-alg\`ebre
commutative $A'$, on a 
$$E \otimes_{A}A' \in [\widehat{T}_{A'}]^{[a,b]}.$$
Par d\'efinition, un objet $E \in \mathcal{M}^{[a,b]}_{T}(A)$ sera dit
\emph{d'amplitude contenu dans l'intervalle $[a,b]$}, ou encore
\emph{d'amplitude $[a,b]$}.

Le lemme suivant montre que 
la condition d'\^etre d'amplitude dans l'intervalle $[a,b]$ 
est une condition locale pour la topologie \'etale. 

\begin{lem}\label{l2}
Le pr\'echamp $\mathcal{M}^{[a,b]}_{T}$ est un champ 
pour la topologie \'etale.
\end{lem}

\textit{Preuve:} 
On commence par remarquer le fait suivant. Soit $A \longrightarrow B$ un morphisme
plat de $k$-alg\`ebres commutatives de type fini. Alors le
dg-foncteur de changement de base
$$\widehat{T}_{A} \longrightarrow \widehat{T}_{B}$$
est t-exact. En effet, sans hypoth\`ese de platitude et simplement par d\'efinition, 
ce dg-foncteur est t-exact \`a droite. Par ailleurs, pour $x \in \widehat{T}_{A}$, 
si on note $p : \widehat{T}_{B} \longrightarrow \widehat{T}_{A}$ le dg-foncteur d'oubli, 
on trouve $p(x\otimes_{A}B) \simeq x \otimes_{A}B$. Comme $B$ est un $A$-module
plat, il est colimite filtrante de $A$-modules projectifs, et ainsi $p(x\otimes_{A}B)$
est une colimite filtrante d'objets de la formes $x \otimes_A M$ avec 
$M$ un $A$-module projectif. Mais un tel $x\otimes_A M$ est un r\'etracte
d'une somme de $x$. Ainsi, on voit que $p(x\otimes_{A}B)$ est dans 
la sous-dg-cat\'egorie pleine de $\widehat{T}_{A}$ engendr\'ee par colimites filtrantes par $x$.
Ainsi, si $x \in \widehat{T}_{A}^{\geq 0}$ il en est de m\^eme de $p(x\otimes_{A}B)$.

Soit maintenant $A \longrightarrow B$ un morphisme fid\`element plat 
entre $k$-alg\`ebres commutatives de type fini. On note $B^{*}$
la $k$-alg\`ebre commutative co-simpliciale co-nerf du morphisme $A \longrightarrow B$, 
de sorte que $B^n=B^{\otimes_{A} n+1}$. On dipose d'une co-augmentation
$A \longrightarrow B^*$. Il faut montrer que le morphisme naturel
$$\mathcal{M}^{[a,b]}_{T}(A) \longrightarrow \lim_{n\in\Delta} \mathcal{M}^{[a,b]}_{T}(B^n)$$
est une \'equivalence. Comme $\mathcal{M}^{[a,b]}_{T}$ est un sous-pr\'echamp
plein du champ $\mathcal{M}_{T}$, ce morphisme est \'equivalent \`a une r\'eunion
de composantes connexes. Il reste \`a montrer que si un objet t-plat
$x \in T_{A}$ est tel que $x\otimes_A B \in T_{B}^{[a,b]}$ alors
$x \in T_{A}^{[a,b]}$. Mais comme $A \mapsto T_{A}$ est un champ
en dg-cat\'egories (voir \cite{to4}), on a une \'equivalence dans $T_{A}$
$$x \simeq \lim_{n\in \Delta} x\otimes_{A}B^{n}.$$
Chaque $B^n$ est un $A$-module plat, et donc colimite filtrante
de $A$-module projectifs. Ainsi, chaque $x\otimes_{A}B^{n}$ est un facteur
d'une somme de $x$ et est donc dans $\widehat{T}_{A}^{[a,b]}$.
Comme $\widehat{T}_{A}^{\geq a}$ est stable par limites dans 
$\widehat{T}_{A}$ il s'en suit que $x\in T_{A}^{\geq a}$.

Par ailleurs, pour montrer que $x \in T_{A}^{\leq b}$, il faut montrer que
pour tout $y\in \widehat{T}_{A}^{>b}$ le complexe de $A$-modules
$\widehat{T}_{A}(x,y)$ est cohomologiquement concentr\'e en degr\'es strictement positifs. 
Comme $A \longrightarrow B$ est fid\`element plat, il suffit de montrer que
$\widehat{T}_{A}(x,y)\otimes_{A}B$ est cohomologiquement concentr\'e en degr\'es 
strictement positifs. 
Or $x$ est compact, et on a donc un quasi-isomorphisme naturel
$$\widehat{T}_{A}(x,y)\otimes_A B \simeq \widehat{T}_{B}(x\otimes_A B,y\otimes_A B).$$
Comme nous avons vu que $\otimes_{A}B$ est t-exact, on a 
$y\otimes_A B \in \widehat{T}_{B}^{>b}$, et donc $H^{i}(\widehat{T}_{B}(x\otimes_A B,y\otimes_A B))=0$
pour $i\leq 0$. Ceci implique que $x \in \widehat{T}_{A}^{\leq b}$ et termine la preuve
du lemme.
\hfill $\Box$ \\

\begin{df}\label{d3}
Soit $T$ une dg-cat\'egorie satur\'ee. Nous dirons qu'une t-structure
sur $T$ est \emph{ouverte} si elle v\'erifie la condition suivante:
pour toute paire d'entiers $a\leq b$, le sous-champ
$$\mathcal{M}_{T}^{[a,b]} \subset \mathcal{M}_{T}$$
est repr\'esentable par une immersion ouverte.
\end{df}

On peut interpr\'eter la d\'efinition ci-dessus comme une tentative 
pour s'assurer d'une forme de semi-continuit\'e des 
foncteurs $H^{i}_t$, analogue \`a la semi-continuit\'e de la dimension
des groupes de cohomologie des complexes parfaits de Grothendieck. \\

Si une dg-cat\'egorie satur\'ee $T$ est munie d'une t-structure ouverte, alors
pour tout $a\leq b$, les champs $\mathcal{M}_{T}^{[a,b]}$ sont 
localement alg\'ebriques et localement de pr\'esentation finie. Par ailleurs, 
si $n=b-a+1$, on voit facilement que $\mathcal{M}_{T}^{[a,b]}$ est 
un $n$-champ d'Artin localement de pr\'esentation finie sur $k$. Ainsi, 
pour $a=b=0$ on trouve un $1$-champ d'Artin $\mathcal{M}_{T}^{[0,0]}$ qui 
classifie les objets compacts du c\oe ur de $T$. Ce champ sera not\'e 
$\mathcal{M}_{T}^{\mathcal{H}}$, et est un champ de modules
d'objets dans la cat\'egorie ab\'elienne $\widehat{\mathcal{H}}$. 
Il peut \^etre d\'ecrit de la mani\`ere suivante. Pour $A$ une $k$-alg\`ebre
commutative notons $\widehat{\mathcal{H}}_A$ la cat\'egorie ab\'elienne
des $A$-modules dans la cat\'egorie $k$-lin\'eaire $\widehat{\mathcal{H}}$. 
L'association $A \mapsto \widehat{\mathcal{H}}_A$ d\'efinit un champ 
en cat\'egories sur le site des $k$-sch\'emas affines muni de la topologie \'etale.
Le champ $\mathcal{M}_{T}^{\mathcal{H}}$ en est le sous-champ en groupo\"\i des
form\'e des objets $E \in \widehat{\mathcal{H}}_A$ qui v\'erifient les deux
conditions suivantes.
\begin{enumerate}

\item L'objet $E$ est \emph{t-plat}: pour tout morphisme de $k$-alg\`ebres
de type fini $A \longrightarrow A'$, l'objet $E\otimes_{A}A' \in \widehat{T}_{A'}$
est dans le c\oe ur de la t-structure induite sur $\widehat{T}_{A'}$.

\item L'objet $E$ est compact dans $[\widehat{T}_A]$. 

\end{enumerate}

Notons que 
le champ $\mathcal{M}_{T}^{\mathcal{H}}$ d\'epend \`a priori de strictement plus que
de la cat\'egorie ab\'elienne, car la condition de compacit\'e pr\'ec\'edente
fait intervenir le plongement $\widehat{\mathcal{H}}_A \hookrightarrow [\widehat{T}_A]$
qui d\'epend \`a priori du choix de $T$. \\

Venons-en \`a notre seconde notion, celle de t-structure parfaite. 

\begin{df}\label{d3'}
Soit $T$ une dg-cat\'egorie triangul\'ee munie d'une t-structure. Nous dirons que la
t-structure est \emph{parfaite} si pour toute $k$-alg\`ebre $A$ r\'eguli\`ere les deux conditions suivantes sont satisfaites.
\begin{enumerate}
\item La t-structure induite sur la dg-cat\'egorie triangul\'ee $T_A=T\widehat{\otimes}_{k}{A}$ 
par extension des scalaires est une t-structure compacte.
\item Le c\oe ur $\widehat{\mathcal{H}}_A$ de la t-structure induite 
est une cat\'egorie ab\'elienne noeth\'erienne: tout sous-objet d'un objet 
$\omega$-petit est $\omega$-petit. 
\end{enumerate}
\end{df}

\section{Espaces Quot}

Nous avons vu la notion de t-structure ouverte sur une dg-cat\'egorie
satur\'ee $T$. Nous allons voir dans ce paragraphe comment cette notion
permet de d\'efinir des g\'en\'eralisation des sch\'emas Quot de Grothendieck
pour des objets de $\mathcal{H}$. Ces espaces Quot seront utilis\'es de mani\`ere
essentielle dans la suite de ce travail pour aboutir \`a notre th\'eor\`eme principal,
mais cette notion poss\`ede aussi un int\'er\^et en soi. 

Soit donc $T$ une dg-cat\'egorie satur\'ee sur k munie d'une t-structure
parfaite et ouverte  (voir les d\'efinitions \ref{d3} et \ref{d3'}). 
Soit $A$ une $k$-alg\`ebre commutative de type fini. 
On note $\mathcal{H}_A$ l'intersection du c\oe ur de la t-structure 
induite sur $[\widehat{T}_A]$ avec $[T_{A}]$ la sous-cat\'egories des objets compacts.
Rappelons qu'un objet $E$ de $\widehat{\mathcal{H}}_A$ est dans 
$\mathcal{H}_A$ si et seulement si pour tout objet compact $K \in T$, le complexe
$T(K,E)$ est un complexe parfait de $A$-modules (prop. \ref{p0} et le fait
que $T$ soit satur\'ee). On suppose que 
$E$ est t-plat: pour tout morphisme de $k$-alg\`ebres commutatives de type fini $A \longrightarrow 
A'$ l'objet $E \otimes_{A}A' \in [\widehat{T}_{A'}]$ est dans le c\oe ur
de la t-structure induite sur $[\widehat{T}_{A'}]$.

Comme dans \cite{tova1}, on notera $\mathcal{M}_{T}^{(1)}$ le champ 
des morphismes dans $T$, classifiant les morphismes de $T^{op}$-dg-modules parfaits. 
Il vient avec deux projections
$$\xymatrix{\mathcal{M}_{T} & \mathcal{M}_{T}^{(1)} \ar[r]^-{t} \ar[l]_-{s} & \mathcal{M}_{T},}$$
de \emph{source} et \emph{but}. La condition de platitude impos\'ee sur $E$ 
d\'efinit un morphisme de champs 
$$E : Spec\, A \longrightarrow \mathcal{M}_{T}^{\mathcal{H}} \subset \mathcal{M}_{T}.$$
On consid\`ere le champ $\mathcal{M}_{T}^{/E}$, d\'efini par le produit fibr\'e
ci-dessous
$$\xymatrix{
\mathcal{M}_{T}^{/E} \ar[r] \ar[d] & \mathcal{M}_{T}^{(1)} \ar[d]^-{t} \\
Spec\, A \ar[r]_-{E} & \mathcal{M}_{T}.}$$
Le champ $\mathcal{M}_{T}^{/E}$ classifie les 
$T^{op}$-dg-modules compacts munis d'un morphisme vers $E$. Sa pr\'esentation par le
produit fibr\'e ci-dessus montre qu'il s'agit d'un champ localement alg\'ebrique 
et localement de type fini sur $k$. 
Le morphisme source d\'efinit
un morphisme induit
$$s : \mathcal{M}_{T}^{/E} \longrightarrow \mathcal{M}_{T}.$$
Nous noterons $\mathcal{M}_{T}^{\mathcal{H}/E}$ le champ d\'efini par le produit fibr\'e
$$\xymatrix{
\mathcal{M}_{T}^{\mathcal{H}/E} \ar[r] \ar[d] & \mathcal{M}_{T}^{\mathcal{H}} \ar[d] \\
\mathcal{M}_{T}^{/E} \ar[r]_-{s} & \mathcal{M}_{T}.}$$
Le champ $\mathcal{M}_{T}^{\mathcal{H}/E}$ est un 1-champ au-dessus de $Spec\, A$, 
dont les sections au-dessus d'une $A$-alg\`ebre commutative de type finie $A'$ 
est le groupo\"\i de des couples $(E',u)$, avec $E' \in \mathcal{H}_{A'}$ 
un objet t-plat, et $u$ un morphisme $u : E' \longrightarrow E\otimes_{A}A'$
dans $\mathcal{H}_{A'}$. En utilisant l'alg\'ebricit\'e de $\mathcal{M}_{T}^{(1)}$
et l'ouverture de la t-structure, il est facile de voir que 
$\mathcal{M}_{T}^{\mathcal{H}/E}$ est un $1$-champ d'Artin, qui se r\'ealise
comme un sous-champ ouvert du champ localement alg\'ebrique $\mathcal{M}_{T}^{/E}$.

Enfin, soit $Quot(E)$ le sous-champ de $\mathcal{M}_{T}^{\mathcal{H}/E}$
form\'e des objets $(E',u)$ comme ci-dessus v\'erifiant la condition suivante: 
le morphisme $u$ est un monomorphisme dans la cat\'egorie ab\'elienne 
$\widehat{\mathcal{H}}_{A'}$, et son conoyau est un objet $Coker(u) \in \widehat{\mathcal{H}}_{A'}$
qui est t-plat. En d'autre termes, $Quot(E)$ est le produit fibr\'e suivant
$$\xymatrix{
Quot(E) \ar[r] \ar[d] & \mathcal{M}_{T}^{\mathcal{H}/E} \ar[d]^-{c} \\
\mathcal{M}_{T}^{\mathcal{H}} \ar[r] & \mathcal{M}_{T},}$$
o\`u $c$ est le morphisme qui envoie une paire $(E',u)$ comme ci-dessus sur le c\^one
de $u$. 

Finalement, on fixe un g\'en\'erateur compact $K$ de $\widehat{T}$, et 
une fonction $\nu : \mathbb{Z} \longrightarrow \mathbb{N}$
\`a support fini. On note $\mathcal{M}_{T}^{\nu}\subset \mathcal{M}_{T}$
le sous-champ ouvert de type fini des $T^{op}$-dg-modules compacts
$F$ tels que (voir \cite[\S 3.3]{tova1})
$$dim_{k}H^{i}(T(K,E))\leq \nu(i) \; \forall \; i\in \mathbb{Z}.$$
Nous noterons de m\^eme $Quot^{\nu}(E)$ le sous-champ ouvert
de $Quot(E)$ form\'e des paires $(E',u)$ dont le c\^one 
de $u$ est dans $\mathcal{M}_{T}^{\nu}$. En d'autre termes, $Quot^{\nu}(E)$
est le produit fibr\'e
$$\xymatrix{
Quot^{\nu}(E) \ar[r] \ar[d] & \mathcal{M}_{T}^{\mathcal{H}/E} \ar[d]^-{c} \\
\mathcal{M}_{T}^{\mathcal{H}}\cap \mathcal{M}_{T}^{\nu} \ar[r] & \mathcal{M}_{T}.}$$

\begin{prop}\label{p1}
Avec les hypoth\`eses, notations et d\'efinitions ci-dessus, les deux assertions suivantes
sont satisfaites.
\begin{enumerate}
\item Le champ 
$Quot^{\nu}(E)$ est repr\'esentable par un espace alg\'ebrique de type fini et s\'epar\'e
sur $Spec\, A$.
\item La projection $Quot(E) \longrightarrow S=Spec\, A$ v\'erifie le crit\`ere
valuatif de propret\'e: pour toute $k$-alg\`ebre commutative $R$ int\`egre lisse et 
de dimension $1$, de corps des fractions $K$, 
le morphisme
$$Quot(E)(R) \longrightarrow Quot(E)(K)\times_{S(K)}S(R)$$
est une \'equivalence.
\end{enumerate} 
\end{prop}

\textit{Preuve:} $(1)$ Les r\'esultats de \cite[\S 3.3]{tova1} impliquent ais\'ement que 
$Quot^{\nu}(E)$ est un espace alg\'ebrique de type fini au-dessus de $Spec\, A$. 
Comme $Quot^{\nu}(E)$ est un sous-champ ouvert de
$Quot(E)$, sa s\'eparation est une cons\'equence du crit\`ere valuatif 
de propret\'e
que nous allons d\'emontrer dans $(2)$.

$(2)$ On fibre le morphisme $Quot(E)(R) \longrightarrow Quot(E)(K)\times_{S(K)}S(R)$
au-dessus de $S(R)$, et on montre que qu'il est bijectif fibre \`a fibre. Cela 
implique que l'on peut remplacer $A$ par $R$. 

On dispose donc de $(E',u)$ un objet de $Quot(E)(R)$, repr\'esent\'e comme une suite
exacte d'objets t-plats 
dans $\mathcal{H}_{R}$
$$\xymatrix{
0 \ar[r] & L \ar[r] & E \ar[r]^-{u} & E' \ar[r] & 0.}$$
On utilise l'adjonction de cat\'egories ab\'eliennes
$$-\otimes_{R}K : \widehat{\mathcal{H}}_{R} \leftrightarrows \widehat{\mathcal{H}}_{K}$$
induite par changement de base et son adjoint \`a droite le foncteur d'oubli.
L'unit\'e de l'ajonction induit un morphisme de suites exactes
dans $\widehat{\mathcal{H}}_{R}$
$$\xymatrix{
0 \ar[r] & L\ar[d]  \ar[r] & E \ar[r]^-{u}\ar[d]  & E' \ar[r] \ar[d] & 0 \\
0 \ar[r] & L\otimes_{R}K \ar[r] & E\otimes_{R}K \ar[r]^-{u_K} & E'\otimes_{R}K \ar[r] & 0.
}$$
On commence par remarquer que les morphismes verticaux sont des monomorphismes. 
En effet, pour tout objet t-plat $E \in \widehat{\mathcal{H}}_{R}$, 
l'unit\'e $v : E \longrightarrow E\otimes_{R}K$ est un monomorphisme. Pour voir cela, 
on consid\`ere le c\^one de $v$ dans $\widehat{T}_{R}$, qui s'exprime comme  
$E\otimes_{R}K/R \simeq \colim_{f\in R-0} E/f$. Ici la colimite 
est prise sur l'ensemble filtrant des \'el\'ements non-nuls de $R$ (ordonn\'e
par la relation de divisibilit\'e), et $E/f$ d\'esigne le c\^one 
de la multiplication $\times f : E \rightarrow E$. Ce c\^one s'exprime aussi 
comme $E\otimes_{R}R/(f)$, et comme $E$ est t-plat ce c\^one est dans
$\mathcal{H}_{R}$. On en d\'eduit que $\colim_{f\in R-0} E/f$, et donc
$E\otimes_{R}K/R$ est dans $\widehat{\mathcal{H}}_{R}$, et ainsi que
le morphisme $v$ est un monomorphisme. 

Le fait que les morphismes verticaux du diagramme pr\'ec\'edent soient 
des monomorphismes implique que le carr\'e
$$\xymatrix{
L \ar[r] \ar[d] & E \ar[d] \\
L\otimes_{R}K \ar[r] & E\otimes_{R}K}$$
est cart\'esien dans la cat\'egorie ab\'elienne $\widehat{\mathcal{H}}_{R}$
(attention, il ne l'est pas dans $\widehat{T}_{R}$). Cela implique clairement 
que $(E,u')$ est l'unique rel\`evement de $(E'\otimes_{R}K,u\otimes_{R}K)$
de $K$ \`a $R$. Ainsi, le morphisme
$Quot(E)(R) \longrightarrow Quot(E)(K)\times_{S(K)}S(R)$ est injectif.

Soit maintenant $(E'_K,u_K)$ un objet de $Quot(E)(K)$, repr\'esent\'e par une suite
exacte dans $\mathcal{H}_{K}$
$$\xymatrix{
0 \ar[r] & L_{K} \ar[r] & E_{K} \ar[r]^-{u_K} & E_{K}' \ar[r] & 0.}$$
On d\'efinit un sous-objet $L\subset E$ dans $\widehat{\mathcal{H}}_{R}$ par le carr\'e
cart\'esien (dans $\widehat{\mathcal{H}}_{R}$) suivant
$$\xymatrix{
L \ar[r] \ar[d] & E \ar[d] \\
L_{K} \ar[r] & E\otimes_{R}K.}$$
Comme la t-structure est noeth\'erienne $L$ est un objet 
$\omega$-petit et donc est dans $\mathcal{H}_{R}$. Par ailleurs, 
le quotient $E/L$ dans $\mathcal{H}_{R}$ est par construction un sous-objet 
de $E'\otimes_{K}R \in \widehat{\mathcal{H}}_{R}$.  Ceci implique 
que $E/L$ est un objet sans torsion: pour tout $f\in R-0$, le c\^one
de la multiplication $\times f : E \longrightarrow E$ 
est dans $\mathcal{H}_{R}$. Ceci implique ais\'ement que 
pour tout $R$-module $M$, l'objet $E\otimes_{R}M$ \`a priori 
dans $\widehat{T}_{R}$, reste dans $\widehat{\mathcal{H}}_{R}$. 
L'objet $E/L$ est donc t-plat, et l'objet $(E/L,u)$, o\`u 
$u : E \longrightarrow E/L$ est la projection naturelle 
est un rel\`evement de $(E'_K,u_K)$ \`a un point de $Quot(E)(R)$. 
\hfill $\Box$ \\

L'espace alg\'ebrique $Quot(E)$ contient deux copies du sch\'ema $Spec\, A$, qui 
correspondent au couple $(E',u)$ avec soit $u$ un isomorphisme
soit $E\simeq 0$. Cela d\'efinit deux morphismes
$$a,b : Spec\,A \rightrightarrows Quot(E),$$ 
de sorte que la projection $Quot(E) \longrightarrow Spec\, A$
en soit un r\'etracte. En particulier, comme $Quot(E)$ est 
s\'epar\'e d'apr\`es la proposition \ref{p1}, les morphismes
$a$ et $b$ sont des immersions ferm\'ees. Ces deux morphismes sont 
aussi clairement des immersions ouvertes, car \^etre nul 
est une condition ouverte pour les complexes parfaits 
(le morphisme $Spec\, k \longrightarrow \mathcal{M}_{T}$ 
correspondant \`a l'objet nul est une immersion ouverte). On trouve 
ainsi une d\'ecomposition canonique d'espaces alg\'ebriques au-dessus de 
$Spec\, A$
$$Quot(E) \simeq   Spec\, A \coprod Quot^{\sharp}(E) \coprod Spec\, A,$$
avec $Quot^{\sharp}(E)$ l'ouvert des paires $(E',u)$
qui sont telles que $u$ n'est pas un isomorphisme et $E'$ n'est pas nul. \\

La d\'ecomposition pr\'ec\'edente et la proposition \ref{p1} implique 
le corollaire suivant.

\begin{cor}\label{cp1}
Avec les notations pr\'ec\'edentes, $Quot^{\sharp}(E)$ est un 
espace alg\'ebrique s\'epar\'e localement de pr\'esentation finie
sur $Spec\, A$, et la projection
$$Quot^{\sharp}(E) \longrightarrow Spec\, A$$
v\'erifie de plus le crit\`ere valuatif de propret\'e. 
\end{cor}

\section{Objets ponctuels et co-engendrement}

On fixe $T$ une dg-cat\'egorie satur\'ee sur $k$. Rappelons (voir section \S 1)
que sur $T$ on dispose de l'auto-\'equivalence de Serre
$$S_T : T \simeq T,$$
telle qu'il existe des quasi-isomorphismes fonctoriels
$$T(x,y)^{\vee} \simeq T(y,S(x)),$$
pout toute paire d'objets $(x,y)$ dans $T$. 

\begin{df}\label{d4}
Un \emph{objet ponctuel de dimension $d \in \mathbb{N}$ dans 
$T$} est un objet $x \in T$ v\'erifiant les conditions suivantes.
\begin{enumerate}

\item Pour tout $i< 0$ on a $H^i(T(x,x)) \simeq 0$. 

\item La dg-alg\`ebre $T(x,x)$ est quasi-isomorphe \`a $Sym_{k}(H^1(T(x,x))[-1])$. 

\item $H^{1}(T(x,x))$ est un $k$-espace vectoriel de dimension $d$. 

\item On a $S_T(x)\simeq x[d]$.

\end{enumerate}
\end{df}

\textbf{Exemples.}
Deux exemples standards d'objets ponctuels. Si $X$ est un espace alg\'ebrique propre et lisse
de dimension $d$, le gratte-ciel en un point $k(x)$, vu comme objet de 
$L_{parf}(X)$ est un objet ponctuel de dimension $d$. Si $A$ est une vari\'et\'e
ab\'elienne de dimension $d$, alors tout fibr\'e en droites de degr\'e $0$ sur
$A$ est un objet ponctuel de dimension $d$.

\begin{df}\label{d5}
Un \emph{syst\`eme de points de dimension $d \in \mathbb{N}$ dans 
$T$} est la donn\'ee d'un ensemble de classes d'\'equivalence d'objets $\mathcal{P}$ de $T$
satisfaisant les conditions suivantes. 
\begin{enumerate}

\item Pour tout $x\in \mathcal{P}$, $x$ est un objet ponctuel de dimension $d$
dans $T$.

\item Pour tout $x,y \in \mathcal{P}$, avec $x\neq y$, on a $T(x,y)\simeq 0$. 

\item Un objet $E \in T$ est nul si et seulement si pour tout $x \in \mathcal{P}$
on a $T(E,x)\simeq 0$. 

\end{enumerate}
\end{df}

Un commentaire sur le point $(3)$ de la d\'efinition pr\'ec\'edente. 
Par dualit\'e, $T(E,x)\simeq 0$ est \'equivalent \`a $T(E,x)^{\vee}\simeq T(S_{T}^{-1}(x),E)\simeq 
0$. Or $T(S_T^{-1}(x),E)\simeq T(x,E)[-d]$ car $x$ est ponctuel de dimension $d$. 
Ainsi, la condition $(3)$ peut aussi s'exprimer par: $E$ est nul 
si et seulement si $T(x,E) \simeq 0$. Cependant, les objets $x \in \mathcal{P}$ ne
forment pas une famille de g\'en\'erateurs compacts de $\widehat{T}$, car il n'est pas
vrai en g\'en\'eral que $T(x,E) \simeq 0$ pour tout $x\in \mathcal{P}$ implique
$E \simeq 0$ pour un objet quelconque de $\widehat{T}$.  \\

\textbf{Exemples.} Pour revenir aux deux exemples pr\'ec\'edents
l'ensemble des gratte-ciel sur $X$ forme un syst\`eme de points de dimension $d$, 
et de m\^eme l'ensemble des fibr\'es en droites de degr\'e z\'ero sur $A$. \\

Bien que les objets consituant un syst\`eme de points dans une
dg-cat\'egorie satur\'ee $T$ ne forment pas 
des g\'en\'erateurs compacts, ils permettent de d\'etecter
les \'equivalences. 
Nous retiendrons le r\'esultat bien connu suivant (voir par
exemple \cite[Lem. 2.15]{or}).

\begin{lem}\label{l3}
Soit $f : T \longrightarrow T'$ un dg-foncteur entre deux 
dg-cat\'egories satur\'ees. Soient $\mathcal{P}$ et 
$\mathcal{P}'$ deux syst\`emes de points de m\^eme dimension $d$
dans $T$ et $T'$. On suppose que les deux conditions suivantes
sont satisfaites.
\begin{enumerate}
\item Le dg-foncteur $f$ envoie $\mathcal{P}$ dans 
$\mathcal{P}'$ et de mani\`ere surjective.
\item Le dg-foncteur $f$ restreint aux objets de $\mathcal{P}$ est 
pleinement fid\`ele: pour tout $x,y \in \mathcal{P}$, le morphisme
$$T(x,y) \longrightarrow T(f(x),f(y))$$
est un quasi-isomorphisme. 
\end{enumerate} 
Alors $f$ est une quasi-\'equivalence. 
\end{lem}

\textit{Preuve:} On note $f_*$ et $f_!$ les adjoints \`a droite et 
\`a gauche de $f : T \longrightarrow T'$ (voir \cite{to4}). Il faut montrer que 
pour tout $z\in T$ et $z' \in T'$ les co-unit\'es d'adjonction
$$f_!f(z) \longrightarrow z \qquad  ff_*(z') \longrightarrow z'$$
sont des \'equivalences dans $T$ et $T'$. Comme $\mathcal{P}$
et $\mathcal{P}'$ 
sont des syst\`emes de points dans $T$ et $T'$, il suffit de voir
que pour tout $x\in \mathcal{P}$ et tout $x' \in \mathcal{P}$
les morphismes induits
$$T'(f_!f(z),x) \longrightarrow T(z,x) \qquad T(x',ff_*(z')) \longrightarrow T'(x',z')$$
sont des quasi-isomorphismes. Par adjonction, il faut v\'erifier que les unit\'es
d'adjonctions
$$x \longrightarrow f_*f(x) \qquad  x' \longrightarrow ff_!(x')$$
sont des \'equivalences. Mais cela est induit par le fait que 
le dg-foncteur $f$ induit une quasi-\'equivalence de la
la sous-dg-cat\'egorie pleine des objets de $\mathcal{P}$ dans $T$ 
vers celle des objets de $\mathcal{P}'$ dans $T'$.
 \hfill $\Box$ \\

Les notions de syst\`eme de points et de t-structure entrent en interaction
\`a travers la notion de co-engendrement, r\'esum\'ee dans la d\'efinition suivante.
Rappelons que si une $t$-structure sur $T$ est parfaite elle induit une
$t$-structure sur $[T]$, qui dispose d'un c\oe ur $\mathcal{H}$. Ce c\oe ur
est une sous-cat\'egorie pleine de $\widehat{\mathcal{H}}$ stable par
sommes finies, r\'etractes, noyaux, conoyaux et extensions. 

\begin{df}\label{d6}
Soit $T$ une dg-cat\'egorie satur\'ee munie d'une t-structure 
parfaite et d'un syst\`eme de points $\mathcal{P}$ de dimension $d$. 

\begin{enumerate}

\item Nous dirons que \emph{le syst\`eme $\mathcal{P}$ co-engendre la t-structure}
si l'on a:
$$(E\in T^{\leq 0}) \Longleftrightarrow (H^{i}(T(E,x))\simeq 0 \; \forall \; x\in \mathcal{P}, \; 
\forall i<0).$$

\item Nous dirons que \emph{le syst\`eme $\mathcal{P}$ co-engendre fortement la t-structure}
s'il co-engendre la t-structure et si de plus pour tout objet $E$ du c\oe ur $\mathcal{H}$
de la t-structure induite sur $T$, on a:
$$(E\simeq 0) \Longleftrightarrow ([E,x] \simeq 0 \; \forall \; x\in \mathcal{P}).$$
\end{enumerate}

\end{df}

Lorsqu'un syst\`eme de points $\mathcal{P}$ co-engendre une t-structure parfaite, cette derni\`ere
est totalement caract\'eris\'ee par la donn\'ee de $\mathcal{P}$. En effet, 
par d\'efinition $\mathcal{P}$ d\'etermine la partie n\'egative $[T]^{\leq 0}=[T^{\leq 0}]$
de la t-structure induite sur $[T]$. On reconstruit alors $\widehat{T}^{\leq 0} \subset \widehat{T}$
comme \'etant la plus petite sous-dg-cat\'egorie pleine contenant 
les objets de $T^{\leq 0}$ et qui est stable par sommes arbitraires et 
par c\^ones (voir prop. \ref{p0} $(2)$). 

\begin{lem}\label{l4}
Soit $T$ une dg-cat\'egorie satur\'ee munie d'une t-structure parfaite co-engendr\'ee
par un syst\`eme de points $\mathcal{P}$ de dimension $d$. 
\begin{enumerate}
\item Tout objet $x \in \mathcal{P}$ appartient au c\oe ur $\mathcal{H}$ de la t-structure
induite sur $[T]$.

\item Le syst\`eme 
$\mathcal{P}$ engendre fortement la t-structure si et seulement si 
$\mathcal{P}$ coincide avec l'ensemble des classes d'isomorphismes d'objets
simples de $\mathcal{H}$. 
\end{enumerate}
\end{lem}

\textit{Preuve:} $(1)$ Soit $x \in \mathcal{P}$. Pour tout $y\in \mathcal{P}$ on 
a, par d\'efinition d'un syst\`eme de points, $H^{i}(T(x,y))\simeq 0$ 
pour $i<0$. Ceci montre que $x \in [T]^{\leq 0}$. Soit maintenant 
$E \in [T]^{\leq -1}$, c'est \`a dire $E[-1] \in [T]^{\leq 0}$. Comme 
$\mathcal{P}$ co-engendre la t-structure on a $H^{i}(T(E[-1],x))\simeq 0$
pour $i<0$, ou en d'autres termes $H^{i}(T(E,x))\simeq 0$ pour $i\leq 0$.
Comme ceci est vrai pour tout $E  \in [T]^{\leq -1}$ on a 
$x \in T^{\geq 0}$. 

$(2)$ Commen\c{c}ons par supposer que $\mathcal{P}$ co-engendre fortement la t-structure. 
Soit $x \in \mathcal{P}$, et soit 
$x \twoheadrightarrow E$ un \'epimorphisme dans $\mathcal{H}$. Si $E\neq 0$, 
alors il existe $y\in \mathcal{P}$ et un morphisme non nul $E \rightarrow y$. 
On a forc\'ement $x=y$, car le morphisme compos\'ee $x \rightarrow y$ est encore
non-nul. Dans ce cas le morphisme compos\'e $x \rightarrow x$ est un scalaire
non-nul de $k$, et donc est un isomorphisme. Cela montre que $x \twoheadrightarrow E$
poss\`ede une section et donc que $E$ est un facteur direct de $x$. Mais comme
$End(x)\simeq k$ cela implique que $E\simeq 0$. Ainsi tous les \'el\'ements
de $x$ sont des objets simples de $\mathcal{H}$. De plus, si $E$ est un objet
simple de $\mathcal{H}$, il est non-nul et donc poss\`ede un morphisme
non-nul $u : E \rightarrow x$ pour un $x \in \mathcal{P}$. On a d\'ej\`a vu
que $x$ \'etait simple, et ainsi le morphisme $u$ et un isomorphisme. Cela montre
que la classe d'isomorphisme de $E$ appartient \`a $\mathcal{P}$. 

Supposons r\'eciproquement que $\mathcal{P}$ soit exactement 
l'ensemble des classes d'isomorphisme d'objets simples de $\mathcal{H}$. 
Soit $E \in \mathcal{H}$ un objet non-nul. Comme $\mathcal{H}$ est noeth\'erienne, 
il existe un objet simple
$x$ de $\mathcal{H}$, et donc un \'el\'ement de $\mathcal{P}$, avec un 
morphisme non-nul $E \rightarrow x$. Cela montre que $\mathcal{P}$ co-engendre
fortement la t-structure. 
\hfill $\Box$ \\

Le point $(2)$ du lemme pr\'ec\'edent implique que lorsqu'une t-structure parfaite et \`a c\oe ur noeth\'erien est 
co-engendr\'ee fortement par un syst\`eme de points, le syst\`eme
$\mathcal{P}$ est lui-m\^eme det\'ermin\'e par la t-structure comme \'etant l'ensemble des
classes d'isomorphismes d'objets simples de $\mathcal{H}$. Il revient ainsi au 
m\^eme, sous ces conditions,  de se donner la t-structure ou 
de se donner le syst\`eme de points $\mathcal{P}$. \\

Pour terminer, on introduit la notion de syst\`emes de points born\'es. 

\begin{df}\label{d7}
Un syst\`eme de points $\mathcal{P}$ dans une dg-cat\'egorie satur\'ee $T$ 
est \emph{born\'e} s'il existe un sous-champ ouvert de type fini sur $k$
$U \subset \mathcal{M}_{T}$ tel que $\mathcal{P}$ soit contenu dans 
le sous-ensemble
$$\pi_{0}(U(k)) \subset \pi_{0}(\mathcal{M}_{T}(k)).$$
\end{df}

Rappelons d'apr\`es \cite[\S 3.3]{tova1} 
la caract\'erisation suivante des syst\`emes de points born\'es.
Un syst\`eme de point $\mathcal{P}$ est born\'e si et seulement si pour tout 
g\'en\'erateur compact $E$ de $\widehat{T}$ il existe 
une fonction $\nu : \mathbb{Z} \longrightarrow \mathbb{N}$ \`a support fini, 
telle que
$$dim_{k}H^{i}(T(E,x)) \leq \nu(i) \; \forall \; i\in \mathbb{Z}, \; \forall x\in \mathcal{P}.$$

\section{L'espace de modules des objets ponctuels}

Pour cette section on fixe une dg-cat\'egorie satur\'ee $T$. On se fixe un syst\`eme de points 
$\mathcal{P}$ de dimension $d$ et une t-structure sur $T$. On suppose que les assertions
suivantes sont satisfaites.

\begin{itemize}

\item La t-structure est parfaite et ouverte.

\item Le famille de points $\mathcal{P}$ co-engendre fortement la t-structure. 

\item La famille de points $\mathcal{P}$ est born\'ee. 

\end{itemize}

Ces donn\'ees permettent de d\'efinir $\mathcal{M}_{\mathcal{P}}$, 
un sous-pr\'echamp de $\mathcal{M}_{T}$ de la mani\`ere suivante. D'apr\`es la d\'efinition
d'\^etre une t-structure ouverte, 
on dispose d'un sous-champ ouvert $\mathcal{M}_{T}^{\mathcal{H}} \subset 
\mathcal{M}_{T}$ form\'e des objets d'amplitude $0$ dans $T$. Pour une $k$-alg\`ebre
commutative de type finie $A$, on d\'efinit un sous-ensemble simplicial plein 
$\mathcal{M}_{\mathcal{P}}(A) \subset \mathcal{M}_{T}^{\mathcal{H}}(A)$, form\'e
des objet $E \in T_A$ satisfaisant \`a la condition suivante: pour 
tout morphisme d'anneaux $A \longrightarrow k$, 
l'objet 
$E\otimes_{A}k \in T$ est dans $\mathcal{P}$. Par d\'efinition de $\mathcal{M}_{T}^{\mathcal{H}}$, 
on sait que pour tout $A \longrightarrow k$  l'objet 
$E\otimes_{A}k$ vit dans le c\oe ur $\mathcal{H}$. 
D'apr\`es le lemme \ref{l4}, la condition pr\'ec\'edente demande de plus que cet objet soit 
un objet simple de $\mathcal{H}$. 

Nous avons d\'efini $\mathcal{M}_{\mathcal{P}}(A)$ pour une $k$-alg\`ebre
de type finie. Nous \'etendons simplement cette d\'efinition \`a toutes
les $k$-alg\`ebres commutatives $A$ en posant
$$\mathcal{M}_{\mathcal{P}}(A) := \colim_{A_{\alpha}\subset A}\mathcal{M}_{\mathcal{P}}(A_{\alpha}),
$$
o\`u la colimite est prise sur l'ensemble filtrant des sous-$k$-alg\`ebres de type fini 
de $A$.

Le principal r\'esultat de repr\'esentabilit\'e est le suivant.

\begin{thm}\label{t1}
Le sous-pr\'echamp $\mathcal{M}_{\mathcal{P}}$ est repr\'esentable
par un $1$-champ alg\'ebrique localement de type fini sur $k$. 
De plus, $\mathcal{M}_{\mathcal{P}}$ est une $\mathbb{G}_{m}$-gerbe
au-dessus de son espace de modules grossier $M_{\mathcal{P}}$, et 
$M_{\mathcal{P}}$ est un espace alg\'ebrique lisse, s\'epar\'e et de type fini sur $k$.
\end{thm}

\textit{Preuve:} Commen\c{c}ons par la repr\'esentabilit\'e de 
$\mathcal{M}_{\mathcal{P}}$. Le sous-pr\'echamp 
$\mathcal{M}_{\mathcal{P}} \subset \mathcal{M}_{T}^{\mathcal{H}}$
est clairement un sous-champ par d\'efinition. Nous allons montrer que l'inclusion
$\mathcal{M}_{\mathcal{P}} \subset \mathcal{M}_{T}^{\mathcal{H}}$ est repr\'esentable
par une immersion ouverte. Pour cela, soit $A$ une $k$-alg\`ebre
commutative et $E \in \mathcal{M}_{T}^{\mathcal{H}}(A)$. Soit $U$ le sous-ensemble
des points ferm\'es $s$ de $S=Spec\, A$ tels que 
$E_{s}:=E\otimes_{A}k(s)$ soit un objet simple de $\mathcal{H}$. 
Comme $\mathcal{P}$ est l'ensemble des classes d'isomorphisme 
d'objets simples de $\mathcal{H}$ il faut montrer que
l'ensemble $U$ est (l'ensemble des points ferm\'es d') un ouvert Zariski du sch\'ema $S$. 

Notons $Z:=S-U$ l'ensemble compl\'ementaire, c'est \`a dire 
l'ensemble des points ferm\'es $s$ de $S$ tels que $E_{s}$ ne soit pas simple dans
$\mathcal{H}$. Ainsi, l'ensemble $Z$ est l'image (au niveau des points ferm\'es)
du morphisme $Quot^{\sharp}(E) \longrightarrow S$. Par la corollaire \ref{cp1} cette image
est stable par sp\'ecialisation. Par ailleurs, comme 
la syst\`eme de points est born\'e, $Z$ est aussi l'image de la projection 
$Quot^{\sharp,\nu}(E) \longrightarrow S$ pour une fonction $\nu : \mathbb{Z} \longrightarrow 
\mathbb{N}$ \`a support fini et 
$Quot^{\sharp,\nu}(E)=Quot^{\sharp}(E)\cap Quot^{\nu}(E)$ est un ouvert
quasi-compact de $Quot^{\sharp}(E)$. 
D'apr\`es la proposition \ref{p1} cette image est donc constructible. 
Ainsi, 
$Z$ est constructible et stable par sp\'ecialisation, et est donc un ferm\'e de $S$. 

Ceci montre, comme annonc\'e, que l'inclusion 
$\mathcal{M}_{\mathcal{P}} \subset \mathcal{M}_{T}^{\mathcal{H}}$ est repr\'esentable
par une immersion ouverte. Comme $\mathcal{M}_{T}^{\mathcal{H}}$ est une $\mathbb{G}_{m}$-gerbe
sur un espace alg\'ebrique de type fini il en est de m\^eme de $\mathcal{M}_{\mathcal{P}}$.
Il nous reste \`a montrer que $\mathcal{M}_{\mathcal{P}}$ est un champ
lisse, et que son espace de modules grossier est s\'epar\'e.  \\

Soit $x \in \mathcal{M}_{\mathcal{P}}(k)$ un point global correspondant 
\`a un objet simple $x\in \mathcal{P}$. Nous allons ici consid\'erer
la version d\'eriv\'ee $\mathbb{R}\mathcal{M}_{T}$ du champ
$\mathcal{M}_{T}$ (qui est not\'e $\mathcal{M}_{T}$ dans \cite{tova1}). Le sous-champ ouvert
$\mathcal{M}_{\mathcal{P}} \subset \mathcal{M}_{T}$ correspond \`a un sous-champ 
d\'eriv\'e ouvert
$$\mathbb{R}\mathcal{M}_{\mathcal{P}} \subset \mathbb{R}\mathcal{M}_{T}.$$
Le champ d\'eriv\'e $\mathbb{R}\mathcal{M}_{\mathcal{P}}$ est un 1-champ d'Artin
d\'eriv\'e dont le tronqu\'e est le champ $\mathcal{M}_{\mathcal{P}}$
$$\mathcal{M}_{\mathcal{P}} \simeq t_0(\mathbb{R}\mathcal{M}_{\mathcal{P}}).$$
Le point $x$ d\'efinit un point global de $\mathbb{R}\mathcal{M}_{\mathcal{P}}$, 
au-quel on peut prendre le complexe tangent $\mathbb{T}_{x}$ (voir \cite{hagII}). Comme montr\'e dans \cite{tova1} on 
trouve un quasi-isomorphisme
$$\mathbb{T}_{x}[-1] \simeq T(x,x).$$
Par ailleurs, on sait que $\mathbb{T}_{x}[-1]$ se promeut en une dg-alg\`ebre
de Lie sur $k$ (voir \cite[\S 4.3]{to1} et \cite{lu3}), et l'identification devient alors un identification 
de dg-alg\`ebres de Lie o\`u la structure de Lie sur le membre
de droite est induite par le commutateur dans la dg-alg\`ebre $T(x,x)$. 
Par d\'efinition des objets ponctuels, 
la dg-alg\`ebre de Lie $T(x,x)$ est quasi-isomorphe \`a $Sym_{k}(V[-1])$, avec 
crochet et diff\'erentielle nulle, et $V=H^{1}(T(x,x))$. D'apr\`es la correspondence
entre champs d\'eriv\'es formels et dg-alg\`ebres de Lie de \cite{lu3}, le compl\'et\'e
formel de $\mathbb{R}\mathcal{M}_{\mathcal{P}}$ en $x$ s'\'ecrit
comme un produit de champs d\'eriv\'es formels
$$B\hat{\mathbb{G}}_{m} \times \widehat{\mathbb{A}}^{d} \times 
F.$$
Le facteur $F$ correspond \`a la dg-lie ab\'elienne $Sym^{\geq 2}(V[-1])$, concentr\'ee
en degr\'es sup\'erieurs \`a $2$. Ainsi, on a $t_0(F)\simeq Spec\, k$. Ceci montre que 
le compl\'et\'e formel du tronqu\'e $\mathcal{M}_{\mathcal{P}}$ en $x$ est 
\'equivalent \`a $B\hat{\mathbb{G}}_{m} \times \widehat{\mathbb{A}}^{d}$. On voit ainsi 
que $\mathcal{M}_{\mathcal{P}}$  est lisse en $x$, et ce pour tout $x$. \\

Pour terminer la preuve du th\'eor\`eme, il nous reste \`a voir que l'espace
de modules grossier  de $\mathcal{M}_{\mathcal{P}}$ est un espace alg\'ebrique
s\'epar\'e. Notons $M_{\mathcal{P}}$ cet espace de modules grossier. C'est un espace
alg\'ebrique lisse et de pr\'esentation finie sur $k$. Pour montrer que
$M_{\mathcal{P}}$ est s\'epar\'e il suffit de montrer que pour 
toute $k$-alg\`ebre $A$ de valuation discr\`ete, que l'on peut supposer
de plus compl\`ete, le morphisme
$M_{\mathcal{P}}(A) \longrightarrow M_{\mathcal{P}}(K)$ est injectif (o\`u 
$K$ est le corps des fractions de $A$). Soient $x$ et $y$ deux \'el\'ements
de $M_{\mathcal{P}}(A)$ d'images \'egales dans $M_{\mathcal{P}}(K)$. 
Comme $A$ est strictement hens\'elien les points $x$ et $y$ se rel\`event en deux \'el\'ements
$x'$ et $y'$ dans $\mathcal{M}_{\mathcal{P}}(A)$. Dans ce cas, $x'$ et $y'$ correspondent
\`a deux objets $E,F \in \mathcal{H}_A$ qui sont t-plats et dont les chang\'es
de base $E\otimes_{A}K$ et $F\otimes_{A}K$ sont isomorphes dans $\mathcal{H}_{K}$. 
Il faut montrer que $E$ et $F$ sont isomorphes dans $\mathcal{H}_A$. 

Soit $\phi_{K} : E\otimes_{A}K \simeq F\otimes_{A}K$ un isomorphisme 
dans $\mathcal{H}_{K}$. Par adjonction ceci correspond \`a un morphisme 
dans $\widehat{\mathcal{H}}_{A}$
$$\phi : E \longrightarrow F\otimes_{A}K\simeq \colim (\xymatrix{
F \ar[r]^-{\pi} & F \ar[r]^-{\pi} & \dots,})$$
o\`u $\pi \in A$ est une uniformisante. Comme $E$ est de pr\'esentation finie
dans $\widehat{\mathcal{H}}_{A}$ le morphisme ci-dessus se factorise par un des facteurs
de la colimite. En d'autre termes, il existe un morphisme
$\phi_A : E \longrightarrow F$
dans $\mathcal{H}_A$ tel que le morphisme induit par changement de base
\`a $K$ soit $\pi^{n}.\phi_K$ pour un certain $n$. Quitte  \`a remplacer
$\phi_K$ par $\pi^{n}.\phi_K$ on supposera donc que l'isomorphisme $\phi_K$ 
provient par changement de base d'un morphisme dans $\mathcal{H}_A$
$$\phi_A : E \longrightarrow F.$$
 
Comme
$T_{A}(E,F)$ est un complexe parfait de $A$-modules et que $A\simeq \lim_{n}A/\pi^n$ 
est complet, on a 
$$T_A(E,F) \simeq \lim_{n}T_A(E,F)/\pi^n \simeq \lim_{n}T_A(E,F/\pi^n),$$
o\`u $K/\pi^n$ est le c\^one du morphisme $\times \pi^n : K \longrightarrow K$. 
Comme les complexes $T_A(E,F/\pi^n)$ sont tous cohomologiquement 
concentr\'es en degr\'es positifs on trouve un isomorphisme naturel
sur leurs $H^0$
$$Hom_{\mathcal{H}_A}(E,F) \simeq \lim_{n}Hom_{\mathcal{H}_A}(E,F/\pi^n).$$
Le morphisme $\phi_A$ est non nul, il existe donc un entier $n$ tel que 
le morphisme induit $\xymatrix{E \ar[r]^-{\phi_A} & F \ar[r] & F/\pi^n}$
soit aussi non-nul. Soit $n$ le plus petit entier tel que 
le morphisme ci-dessus soit non nul, de sorte que 
l'on ait $\phi_A = \pi^{n-1}.\psi_A$, avec $\psi_A : E \longrightarrow F$
un morphisme tel que le morphisme induit sur le corps
r\'esiduel $\psi_k :=\psi_A \otimes_{A}k : E\otimes_A k \longrightarrow F\otimes_A k$ 
soit un morphisme non nul de $\mathcal{H}$. Comme les objets 
$ E\otimes_A k$ et $F\otimes_A k$ sont suppos\'es simples, le morphisme
$\psi_k$ est un isomorphisme. Enfin, comme $E$ et $F$ sont des 
$T^{op}\otimes_{k}A$ dg-modules compacts, ils sont parfaits
comme $A$-dg-modules. Ainsi, le morphisme $\psi_A : E \longrightarrow F$
est un isomorphisme par Nakayama pour les complexes parfaits. 
Ceci termine de d\'emontrer que $E$ et $F$ sont
isomorphes dans $\mathcal{H}_A$ et donc termine la preuve du th\'eor\`eme
\ref{t1}. \hfill $\Box$ \\

\section{Th\'eor\`eme de reconstruction}

Pour une dg-cat\'egorie satur\'ee $T$, on dispose du champ
$\mathcal{M}_{T}$ des objets de $T$, ainsi que d'un dg-foncteur canonique
$$\phi : T^{op} \longrightarrow L_{parf}(\mathcal{M}_{T}),$$
de la dg-cat\'egorie oppos\'ee \`a celle de $T$ vers celle
des complexes parfaits sur le champs $\mathcal{M}_{T}$. Le morphisme
$\phi$ est par d\'efinition celui correspondant au $T^{op}$-dg-module
universel sur $\mathcal{M}_{T}$. Pour un objet $x \in T$, 
le complexe parfait  $\phi(x)$ sur $\mathcal{M}_{T}$ peut \^etre
d\'ecrit de la mani\`ere suivante. Soit $u : Spec\, A \longrightarrow \mathcal{M}_{T}$
un morphisme correspondant \`a un $T^{op}\otimes_k A$-dg-module parfait $E$, 
que l'on voit comme un dg-foncteur $E : T^{op} \longrightarrow L_{parf}(A)$. 
Par d\'efinition, on dispose d'un identification naturelle de complexes
parfaits de $A$-modules
$$E(x) \simeq u^{*}(\phi(x)).$$
L'identification ci-dessus peut aussi s'\'ecrire 
$$T(x,E) \simeq u^{*}(\phi(x)),$$
o\`u $x \in T$ est vu comme un $T^{op}$-dg-module par le plongement de Yoneda $T \hookrightarrow
\widehat{T}$. En particulier, lorsque $A=k$ et $u : Spec\, k \longrightarrow \mathcal{M}_{T}$
correspdond \`a un objet $y \in T$, on trouve la formule suivante pour la fibre en $u$ 
de $\phi(x)$
$$\phi(x)_{y}:=u^{*}(\phi(x)) \simeq T(x,y) \in L_{parf}(k).$$

On suppose maintenant que l'on est sous les hypoth\`eses du th\'eor\`eme \ref{t1}: on 
se fixe un syst\`eme de points $\mathcal{P}$ de dimension $d$ et une t-structure sur $T$
qui satisfont aux conditions suivantes.

\begin{itemize}

\item La t-structure est parfaite et ouverte.

\item La famille de points $\mathcal{P}$ co-engendre fortement la t-structure. 

\item La famille de points $\mathcal{P}$ est born\'ee. 

\end{itemize}

Par le th\'eor\`eme \ref{t1} on dispose d'une sous-champ ouvert
$j : \mathcal{M}_{\mathcal{P}}\subset \mathcal{M}_{T}$, qui est 
une $\mathbb{G}_{m}$-gerbe sur un espace alg\'ebrique lisse et s\'epar\'e
$M_{\mathcal{P}}$. En composant 
$\phi$ avec la restriction le long de $j$ on trouve un dg-foncteur
$$\phi_{\mathcal{P}} :\xymatrix{T^{op} \ar[r]^-{\phi} & L_{parf}(\mathcal{M}_{T}) \ar[r]^-{j^*} & 
L_{parf}(\mathcal{M}_{\mathcal{P}}).
}$$
En utilisant la description des fibres des objets $\phi(x)$ que l'on donne ci-dessus,
on voit que le dg-foncteur $\phi_{\mathcal{P}}$ se factorise par 
la sous-dg-cat\'egorie $L_{parf}^{\chi=1}(\mathcal{M}_{\mathcal{P}}) \subset 
L_{parf}(\mathcal{M}_{\mathcal{P}})$ des objets de poids $1$.

\begin{df}\label{d8}
Avec les conditions et les notations pr\'ec\'edentes, le \emph{dg-foncteur
de d\'ecomposition associ\'e \`a la paire $(T,\mathcal{P})$} est 
le dg-foncteur
$$\phi_{\mathcal{P}} : T^{op} \longrightarrow L_{parf}^{\chi=1}(\mathcal{M}_{\mathcal{P}})$$
d\'efini ci-dessus.
\end{df}

Intuitivement, l'espace alg\'ebrique $M_{\mathcal{P}}$ est un espace de modules
pour les objets de $\mathcal{P}$, et le dg-foncteur $\phi_{\mathcal{P}}$
d\'ecompose chaque objet $x\in T$ en un complexe parfait $\phi_{\mathcal{P}}(x)$
au-dessus de $M_{\mathcal{P}}$ dont la fibre en $y \in M_{\mathcal{P}}$
est le complexe $T(x,y)$. La non-existence de l'objet universel sur
$M_{\mathcal{P}}$, control\'ee par la gerbe $\mathcal{M}_{\mathcal{P}}$, implique
que $\phi_{\mathcal{P}}(x)$ existe uniquement comme un complexe parfait tordu sur 
$M_{\mathcal{P}}$. D'une certaine fa\c{c}on, le dg-foncteur 
$\phi_{\mathcal{P}}$ doit \^etre consid\'er\'e comme un foncteur de \emph{d\'ecomposition 
spectrale} des objets de $T$ au-dessus de l'espace des modules des objets simples
de la cat\'egorie ab\'elienne $\mathcal{H}$. \\

Nous arrivons enfin au th\'eor\`eme principal de ce travail.

\begin{thm}\label{t2}
Soit $T$ une dg-cat\'egorie satur\'ee munie d'un syst\`eme de points $\mathcal{P}$
v\'erifiant les conditions ci-dessous. 

\begin{itemize}

\item La t-structure est parfaite et ouverte.

\item La famille de points $\mathcal{P}$ co-engendre fortement la t-structure. 

\item La famille de points $\mathcal{P}$ est born\'ee. 

\end{itemize}

Le dg-foncteur de d\'ecomposition 
$$\phi_{\mathcal{P}} : T^{op} \longrightarrow L_{parf}^{\chi=1}(\mathcal{M}_{\mathcal{P}})$$
est conservatif. C'est une quasi-\'equivalence si et seulement si l'espace
alg\'ebrique $M_{\mathcal{P}}$ est propre sur $Spec\, k$. 
\end{thm}

\textit{Preuve:} La conservativit\'e est 
une cons\'equence directe de la formule donnant la fibre de $\phi_{\mathcal{P}}(x)$ en 
un point global $y \in \mathcal{M}_{\mathcal{P}}(k)$
$$\phi_{\mathcal{P}}(x)_{y} \simeq T(x,y),$$
et de la d\'efinition d'un syst\`eme de points $\mathcal{P}$ (voir d\'efinition \ref{d5}).

Supposons pour commencer que $\phi_{\mathcal{P}}$ soit une quasi-\'equivalence.
Comme $T$ est propre, le fait que $\phi_{\mathcal{P}}$ est une \'equivalence
implique que $L_{parf}^{\chi=1}(\mathcal{M}_{\mathcal{P}})$ est une dg-cat\'egorie propre. 

\begin{lem}\label{l5}
Soit $X$ un espace alg\'ebrique de type fini et s\'epar\'e sur $Spec\, k$, 
et $\mathcal{X} \longrightarrow X$ une $\mathbb{G}_{m}$-gerbe sur $X$. 
Si la dg-cat\'egorie $L_{parf}^{\chi=1}(\mathcal{X})$ est propre alors
$X$ est propre sur $Spec\, k$.
\end{lem}

\textit{Preuve du lemme:} Supposons que $X$ ne soit pas propre. Alors, 
il existe une courbe affine lisse $C$ et un morphisme fini
$p : C \longrightarrow X$. 
Par \cite[Cor. 5.2]{to4} (on se ram\`ene ais\'ement au cas o\`u $X$ 
sch\'ema en  utilisant le lemme de Chow, comme par exemple dans 
\cite[Prop. B1]{tova2}). Choississons un g\'en\'erateur compact
$K \in L_{parf}^{\chi=1}(\mathcal{X})$ de la dg-cat\'egorie d\'eriv\'ee
tordue $L^{\chi=1}(\mathcal{X})$, que nous supposerons aussi
g\'en\'erateur compact local. Notons $\mathcal{A}_X:=\mathbb{R}\underline{End}(K) \in 
L_{parf}(X)$ la dg-alg\`ebre d'Azumaya sur $X$ correspondante.

Notons $E:=p^*(\mathcal{A}_X) \in L_{parf}(C)$. Pour tout $V \in L(C)$, on a 
$$\mathbb{R}\Gamma(C,V\otimes E) \simeq \mathbb{R}\Gamma(X,p_*(V)\otimes \mathcal{A}_X) \simeq
\mathbb{R}\underline{Hom}(K,p_*(V)\otimes K).$$
Ainsi, comme $L_{parf}^{\chi}(\mathcal{X})$ est propre, 
le complexe $\mathbb{R}\Gamma(C,V\otimes E)$ est parfait sur $k$. Comme $C$ est affine et que 
$\mathcal{A}_{X}$ est un g\'en\'erateur compact local, $E$ est un g\'en\'erateur
compact de $L(C)$, et le faisceau structural $\mathcal{O}_C$ appartient donc 
\`a la sous-cat\'egorie triangul\'ee \'epaisse engendr\'ee par $E$. Cela implique 
donc que $\mathbb{R}\Gamma(C,V)$ est parfait pour tout complexe parfait $V$ sur $C$, 
et en particulier pour $V=\mathcal{O}_C$.
Ceci est une contradiction car $C$ est une courbe affine. 
\hfill $\Box$ \\

Il nous reste \`a montrer que si $M_{\mathcal{P}}$ est propre alors
$\phi_{\mathcal{P}}$ est une \'equivalence. Tout d'abord on sait que 
$L_{parf}^{\chi=1}(\mathcal{X})$ est une dg-cat\'egorie satur\'ee (voir 
\cite[Cor. 5.4]{to4} $(4)$ pour $X$ un sch\'ema, et le cas des espaces alg\'ebriques
se ram\`ene au cas des sch\'emas par le lemme de Chow comme dans \cite[Prop. B1]{tova2}).
Par ailleurs, il existe un syst\`eme de points $\mathcal{P}'$  de dimension $d$ 
dans $L_{parf}^{\chi}(\mathcal{X})$ d\'efini par les points de $X$ de la mani\`ere
suivante. Pour $x \in X(k)$, on consid\`ere la gerbe r\'esiduelle 
$\mathcal{X}_x:=\mathcal{X}\times_{X}\{x\}\simeq
B\mathbb{G}_m$, ainsi que le morphisme naturel $j_x : \mathcal{X}_x \longrightarrow \mathcal{X}$.
On note $L$ la repr\'esentation de rang $1$ de $\mathbb{G}_m$ et de poids $1$ que 
l'on voit comme un objet de $L[-d]\in L_{parf}^{\chi=1}(\mathcal{X}_x)$. 
Par d\'efinition, $\mathcal{P}'$ est l'ensemble des classes d'\'equivalence 
des objets de la forme $(j_x)_*(L[-d]) \in L_{parf}^{\chi=1}(\mathcal{X})$ (il est donc
en bijection naturelle avec les points de $X$). Il est facile de v\'erifier que
$\mathcal{P}'$ est un syst\`eme de points dans $L_{parf}^{\chi=1}(\mathcal{X})$
au sens de la d\'efinition \ref{d5}. 

Nous allons maintenant appliquer le lemme \ref{l3}
au dg-foncteur de d\'ecomposition
$$\phi_{\mathcal{P}} : T^{op} \longrightarrow L_{parf}^{\chi=1}(\mathcal{X}),$$
avec les syst\`emes de points $\mathcal{P}$ et $\mathcal{P}'$. Il faut 
montrer que $\phi_{\mathcal{P}}$ induit une quasi-\'equivalence
entre les deux sous-dg-cat\'egories pleine de $T^{op}$ et $L_{parf}^{\chi=1}(\mathcal{X})$
d\'efinies par les ensembles de points $\mathcal{P}$ et $\mathcal{P}'$. 

Pour cela, on rappelle que si $x,y\in T^{op}$, alors $\phi_{\mathcal{P}}(x)_y\simeq T(x,y)$, 
o\`u $\phi_{\mathcal{P}}(x)_y=(j_y)^*(\phi_{\mathcal{P}}(x))$ et 
$j_y : \mathcal{X}_y \longrightarrow \mathcal{X}$ est la gerbe r\'esiduelle en $y$. 
Cette identification pour $x=y$ nous donne un morphisme naturel 
dans $L(\mathcal{X}_x)$
$$(j_x)^*(\phi_{\mathcal{P}}(x)) \simeq T(x,x) \longrightarrow H^{d}(T(x,x))[-d]\simeq L[-d].$$
Par adjonction ceci d\'efinit des morphismes dans $L(\mathcal{X})$
$$\phi_{\mathcal{P}}(x) \longrightarrow (j_x)_*(L[-d]),$$
qui sont des \'equivalences. Ainsi, le dg-foncteur $\phi_{\mathcal{P}}$
identifie les ensembles $\mathcal{P}$ et $\mathcal{P}'$, en identifiant un point $x\in \mathcal{P}$
\`a l'objet $(j_x)_*(L)[-d] \in L(\mathcal{X})$. Cette description implique aussi 
que $\phi_{\mathcal{P}}$ est pleinement fid\`ele lorsque il est restreint 
aux objets de $\mathcal{P}$. 
\hfill $\Box$ \\

\section{Deux exemples}

Nous donnons pour terminer deux exemples de situations qui illustrent 
le th\'eor\`eme \ref{t2}. 
Cependant, nous ne donnons pas les d\'etails
des preuves que les hypoth\`eses sont effectivemet v\'erif\'ees dans chacun 
des cas. Il existe par ailleurs de nombreux autres exemples qui m\'eriteraient 
d'\^etre \'etudi\'es, comme par exemple les flops en dimension $3$ de 
\cite{br}. \\

\textbf{Gerbes sur un groupe fini.} Soit $Z$ un espace alg\'ebrique propre et lisse
sur $Spec\, k$. Soit $\mathcal{Z}$ un champ alg\'ebrique
muni d'un morphisme $\pi : \mathcal{Z} \longrightarrow Z$. On suppose que 
localement pour la topologie \'etale sur $Z$ le morphisme $\pi$ est \'equivalent 
\`a la projection $Z \times BG \longrightarrow Z$ pour un groupe fini 
$G$ (qui n'est bien d\'efini que localement). On supposera
pour simplifier que $Z$ est connexe, de sorte que la classe
de conjugaison du groupe $G$ ne d\'epende pas du choix du point $z \in Z$.
On consid\`ere
$T:=L_{parf}(\mathcal{Z})$ la dg-cat\'egorie des complexes parfaits sur 
$\mathcal{Z}$. D'apr\`es \cite[Cor. 5.4]{to4} (et l'astuce
du lemme de Chow de \cite[Prop. B1]{tova2}) on sait que $T$ est une dg-cat\'egorie propre et lisse sur 
$k$. 

Nous d\'efinissons un syst\`eme de points $\mathcal{P}$ de dimension $d=dim\, Z$ 
de la mani\`ere suivante. Soit $z \in Z$ un point ferm\'e, et notons
$\mathcal{Z}_z:=\pi^{1}(z)$ la gerbe r\'esiduelle du champ 
$\mathcal{Z}$ en $z$. Le champ $\mathcal{Z}_z$ est, de mani\`ere non-canonique, 
\'equivalente \`a $BG$ pour $G$ un groupe fini. Cela implique que 
la cat\'egorie des faisceaux quasi-coh\'erents $QCoh(\mathcal{Z}_z)$ s'\'ecrit 
comme un produit fini de cat\'egorie de k-espaces 
vectoriels $\prod_{\rho \in I(z)}Vect(k)$, o\`u $I(z)$ est un ensemble fini
de repr\'esentant d'objets simples dans $QCoh(\mathcal{Z}_z)$. Si l'on choisit 
un \'equivalence $\mathcal{Z}_z \simeq BG$, l'ensemble $I(z)$ peut s'identifier
\`a un ensemble de repr\'esentants des classes d'isomorphisme 
de repr\'esentations lin\'eaires de $G$ sur $k$. Nous noterons

On pose
$$\mathcal{P}:=\coprod_{z \in Z} \coprod_{\rho \in I(z)}(j_{z})^*(\rho),$$
o\`u $j_z : \mathcal{Z}_z \hookrightarrow \mathcal{Z}$
est l'immersion canonique. L'ensemble $\mathcal{P}$ est un ensemble
de faisceaux coh\'erents sur le champ $\mathcal{Z}$ que l'on consid\`ere
comme des objets de $T$. La paire $(T,\mathcal{P})$ v\'erifie les hypoth\`eses
du th\'eor\`eme \ref{t2}. L'espace alg\'ebrique $M_{\mathcal{P}}$ est ici 
un espace alg\'ebrique fini et \'etale sur $Z$, dont la fibre
en $z \in Z$ est l'ensemble fini des classes d'isomorphismes 
de repr\'esentations lin\'eaires du groupe $G$. L'espace $M_{\mathcal{P}}$, ainsi que la classe
$\alpha \in H^2_{et}(Z,\mathbb{G}_m)$ de la $\mathbb{G}_m$-gerbe $\mathcal{M}_{\mathcal{P}} 
\longrightarrow M_{\mathcal{P}}$
peuvent se d\'ecrire explicitement de la mani\`ere suivante. On choisit 
un point de base $z\in Z(k)$, et une identification $\mathcal{Z}_{z}\simeq BG$, 
de sorte que $\mathcal{Z} \longrightarrow Z$ s'identifie \`a une forme tordue, 
pour la topologie \'etale de $BG$. La projection $\pi$ est donc
class\'ee par un morphisme de champs
$$Z \longrightarrow BAut(BG),$$
o\`u $Aut(BG)$ est le 2-groupe des auto-\'equivalences
de $BG$. On se fixe une identification $QCoh(BG) \simeq \oplus_{I}Vect(k)$, avec 
$I$ un ensemble fini qui s'identifie aux classes d'isomorphisme  
de repr\'esentations irr\'ecutibles de $G$ sur $k$. Le champ en groupe
des auto-\'equivalences de $\oplus_{I}Vect(k)$ est un produit
semi-direct $\Sigma_{I} \rtimes B(\mathbb{G}_m)^{I}$. On dipose ainsi d'un morphisme 
de 2-champs $BAut(BG) \longrightarrow B(\Sigma_{I} \rtimes B(\mathbb{G}_m)^{I})$, et donc
d'un morphisme classifiant
$$q : Z \longrightarrow B(\Sigma_{I} \rtimes B(\mathbb{G}_m)^{I}).$$
Ce morphisme d\'etermine le champ en cat\'egories $\mathcal{O}_{Z}$-li\'eaires $\pi_*(QCoh)$ 
sur le petit site \'etale $Z_{et}$. 

On a une suite de fibrations de $2$-champs
$$\xymatrix{B(\mathbb{G}_m)^{I} \ar[r] & B(\Sigma_{I} \rtimes B(\mathbb{G}_m)^{I}) \ar[r]
& B\Sigma_{I}.}$$
La donn\'ee de la projection $p : Z \longrightarrow B(\Sigma_{I})$ d\'etermine un rev\^etement 
\'etale fini $r : Z' \longrightarrow Z$ de fibre $I$. Par ailleurs, 
le rel\`evement de $p$ en un morphisme $q$
$$\xymatrix{Z \ar[r]^-{q} \ar[rd]_-{p} & B(\Sigma_{I} \rtimes B(\mathbb{G}_m)^{I}) \ar[d] \\
 & B(\Sigma_I),}$$
d\'etermine une classe de cohomologie $\alpha \in H^{2}_{et}(Z,A)$, o\`u $A$ est une 
forme tordue du faisceau $\mathbb{G}_m^{I}$. Cette forme tordue n'est autre 
que $r_*(\mathbb{G}_m)$, l'image direct du faisceau $\mathbb{G}_m$ par le
rev\^etement fini $r : Z' \longrightarrow Z$. On trouve donc ainsi une classe
$$\alpha \in H^{2}_{et}(Z',r_*(\mathbb{G}_{m})) \simeq H^{2}_{et}(Z',\mathbb{G}_m).$$

Avec ces notations, l'espace alg\'ebrique $M_{\mathcal{P}}$ s'identifie naturellement \`a
$Z'$, et $\alpha$ ainsi construit est la classe de la gerbe 
$\mathcal{M}_{\mathcal{P}} \longrightarrow M_{\mathcal{P}}$. 

Le th\'eor\`eme \ref{t2} induit ainsi une \'equivalence naturelle
$$L_{parf}(\mathcal{Z}) \simeq L_{parf}^{\alpha}(Z'),$$
o\`u $L_{parf}^{\alpha}(Z')$ est la dg-cat\'egorie des complexes
parfaits sur $Z'$ tordus par la classe $\alpha$. Cet \'enonc\'e est en r\'ealit\'e
vrai sans les hypoth\`eses de lissit\'e et de propret\'e sur $Z$, comme cela peut se
voir de mani\`ere directe. \\

\textbf{Correspondance de McKay.} Soit $X$ une surface projective, lisse et connexe
sur $k$ munie d'une trivialisation $\omega : \omega_X \simeq \mathcal{O}_{X}$. On se
donne un groupe fini $G$ qui op\`ere sur $X$ en fixant la trivialisation 
$\omega$. On consid\`ere $T:=L_{parf}([X/G])$ la dg-cat\'egorie des 
complexes parfaits sur le champ quotient $[X/G]$, qui s'identifie naturellement 
\`a la dg-cat\'egorie des complexes parfaits $G$-\'equivariants sur $X$. 
On d\'efinit un syst\`eme de points $\mathcal{P}$ de dimension 
$2$ dans $T$ de la mani\`ere suivante. On consid\`ere les coh\'erents
$G$-\'equivariants $\mathcal{F}$ sur $X$ qui v\'erifient les deux conditions suivantes.
\begin{itemize}

\item $\mathcal{F}$ est le faisceau structural d'un sous-sch\'ema fini  de $X$
de longueur $\sharp G$.

\item La repr\'esentation $G$ sur $\Gamma(X,\mathcal{F})$ est isomorphe
\`a la repr\'esentation r\'eguli\`ere. 

\end{itemize}

Ces faisceaux coh\'erents forment une famille d'objets dans $T$, et 
nous notons
$\mathcal{P}$ l'ensemble de leurs classes d'\'equivalences
dans $T$.
Le couple 
$(T,\mathcal{P})$ v\'erifie les conditions du th\'eor\`eme \ref{t2}.
Dans cet exemple $M_{\mathcal{P}} \simeq Z$ est propre, et la 
$\mathbb{G}_m$-gerbe est ici triviale. Ainsi, 
le th\'eor\`eme \ref{t2} implique l'existence d'une \'equivalence
$$L_{parf}([X/G]) \simeq L_{parf}(Z),$$
qui est une incarnation de la correspondance de McKay (voir par 
exemple \cite{brki}). \\


\begin{thebibliography}{50}

\bibitem[Ab]{abou} Abouzaid, M. \emph{
Family Floer cohomology and mirror symmetry.} ICM lecture, Seoul 2014. 

\bibitem[BBD]{bbd}  Beilinson, A. A.; Bernstein, J.; Deligne, P. \emph{Faisceaux pervers.}  Analyse et topologie sur les espaces singuliers  I (Luminy, 1981), 5-171, Ast\'erisque, 100, Soc. Math. France, Paris, 1982.

\bibitem[Br]{br} Bridgeland, T. \emph{
Flops and derived categories.}
Invent. Math. \textbf{147} (2002), no. 3, 613-632.

\bibitem[Br-Ki]{brki}
Bridgeland, T.; King, A. \emph{
The McKay correspondence as an equivalence of derived categories.}
J. Amer. Math. Soc. \textbf{14} (2001), no. 3, 535-554 

\bibitem[Ho]{ho} Hovey, M. \emph{Model Categories.} Mathematical Surveys and Monographs, vol 63, Amer. Math. Soc., Providence, 1998.  

\bibitem[Ke]{ke} Keller, B. \emph{On differential graded categories.}
International Congress of Mathematicians. Vol. II, 151-190, Eur. Math. Soc., Z\"urich, 2006.

\bibitem[Ke-Vo]{kv} Keller, B. ; Vossieck, D. \emph{Aisles in derived categories.} 
Bull. Soc. Math. Belg. \textbf{40} (1988), 239-253. 

\bibitem[Lu1]{lu1} Lurie, J. \emph{The "DAG" series.} Available at the author home page
http://www.math.harvard.edu/~lurie/

\bibitem[Lu2]{lu2} Lurie, J. \emph{On the Classification of Topological Field Theories.} Current Developments in Mathematics, vol. 2008, 129-280, Int. Press, Somerville, 2009. 

\bibitem[Lu3]{lu3} Lurie, J. 
\emph{Moduli problems for ring spectra.} Proceedings of the International Congress of 
Mathematicians. Volume II, 1099-1125, Hindustan Book Agency, New Delhi, 2010. 

\bibitem[Lu4]{ha} Lurie, J. \emph{Higher Algebra.} Available at the author home page
http://www.math.harvard.edu/~lurie/

\bibitem[Or]{or}
Orlov, D. \emph{Equivalences of derived categories and K3 surfaces.}
Journal of Mathematical Sciences.
Vol. \textbf{84}, no. 5, 1997.


\bibitem[Pa-To-Va-Ve]{ptvv} 
Pantev, T.; To\"en, B.; Vaqui\'e, M.; Vezzosi, G. \emph{Shifted symplectic structures.}
Publ. Math. Inst. Hautes \'Etudes Sci. \textbf{117} (2013), 271-328.

\bibitem[Ro]{rouq} Rouquier, R. 
\emph{Cat\'egories d\'eriv\'ees et g\'eom\'etrie birationnelle (d'apr\`es Bondal, Orlov, Bridgeland, 
Kawamata et al.)}.  S\'eminaire Bourbaki. Vol. 2004/2005. Ast\'erisque No.\textbf{307} 
(2006), Exp. No. 946, viii, 283-307.

\bibitem[Ta1]{ta} Tabuada, G. \emph{Une structure de cat\'egorie de mod\`eles de Quillen sur la 
cat\'egorie des 
dg-cat\'egories.}
C. R. Math. Acad. Sci. Paris \textbf{340} (2005), no. 1, 15-19. 

\bibitem[Ta2]{ta2} Tabuada, G. \emph{A guided tour through the garden of noncommutative motives.} 
Topics in noncommutative geometry, 259-276, Clay Math. Proc., 16, Amer. Math. Soc., Providence, RI, 
2012.

\bibitem[Ta3]{ta3} Tabuada, G. \emph{Differential graded versus simplicial categories.} 
Topology Appl. \textbf{157} (2010), no. 3, 563-593.

\bibitem[To1]{to1} To\"en, B. 
\emph{Derived Algebraic Geometry.} EMS Surv. Math. Sci. 1 (2014), no. \textbf{2}, 153-245.

\bibitem[To2]{to2} To\"en, B. \emph{The homotopy theory of $dg$-categories and derived Morita 
theory.}  Invent. Math. \textbf{167} (2007), no. 3, 615-667. 

\bibitem[To3]{to3} To\"en, B. \emph{Lectures on $dg$-categories.} Topics in Algebraic and Topological 
*K*-theory. lectures Notes in mathematics, vol. 2008, Springer, 2011. 

\bibitem[To4]{to4} To\"en, B. \emph{Derived Azumaya algebras and generators for twisted derived 
categories.}
 Invent. Math. \textbf{189} (2012), no. 3, 581-652. 

\bibitem[To-Va1]{tova1} To\"en, B.; Vaqui\'e, M. \emph{Moduli of objects in dg-categories.}
 Ann. Sci. \'Ecole Norm. Sup. (4) \textbf{40} (2007), no. 3, 387-444.

\bibitem[To-Va2]{tova2} To\"en, B.; Vaqui\'e, M. \emph{Alg\'ebrisation des vari\'et\'es analytiques 
complexes et cat\'egories d\'eriv\'ees.}
Math. Ann. \textbf{342} (2008), no. 4, 389-831.

\bibitem[To-Ve1]{hagI} To\"en, B.; Vezzosi, G. 
\emph{Homotopical algebraic geometry. I. Topos theory.} Adv. Math. 
\textbf{193} (2005), no. 2, 257-372.

\bibitem[To-Ve2]{hagII} To\"en, B.; Vezzosi, G. 
\emph{Homotopical algebraic geometry. II. Geometric stacks and applications.} 
Mem. Amer. Math. Soc. \textbf{193} (2008), no. 902, x+224 pp.


\end{thebibliography}
\end{document}